\input amstex
\documentstyle{amsppt}
%----------------------------------------------------------------
% Title:     Comparison of two classifications of a class of ODE's
%            in the first case of intermediate degeneration.
% Author:    Ruslan Sharipov
% Comments:  AmSTeX, 26 pages, amsppt style, US letter size, double-sided
% MSC-class: 34A34, 34A26, 34C20, 34C14
%----------------------------------------------------------------
%           Replacement for output macro definition
%
\catcode`@=11
\redefine\output@{%
  \def\break{\penalty-\@M}\let\par\endgraf
  \ifodd\pageno\global\hoffset=105pt\else\global\hoffset=8pt\fi  
  \shipout\vbox{%
    \ifplain@
      \let\makeheadline\relax \let\makefootline\relax
    \else
      \iffirstpage@ \global\firstpage@false
        \let\rightheadline\frheadline
        \let\leftheadline\flheadline
      \else
        \ifrunheads@ %\let\makefootline\relax
        \else \let\makeheadline\relax
        \fi
      \fi
    \fi
    \makeheadline \pagebody \makefootline}%
  \advancepageno \ifnum\outputpenalty>-\@MM\else\dosupereject\fi
}
\def\Beta{\mathchar"0\hexnumber@\rmfam 42}
\catcode`\@=\active
%----------------------------------------------------------------
\nopagenumbers
\chardef\textvolna='176

\font\tencyr=wncyr10

% \font\tenbfsl=cmbxsl10
\chardef\bigalpha='013
\def\negskp{\hskip -2pt}

\def\compos{\,\raise 1pt\hbox{$\sssize\circ$} \,}

\chardef\degree="5E
\def\tr{\operatorname{tr}}
%\def\id{\operatorname{id}}
%\font\eightrm=cmr8
%\def\LT{\operatorname{\text{\eightrm LT}}}
%\def\LM{\operatorname{\text{\eightrm LM}}}
%\def\LC{\operatorname{\text{\eightrm LC}}}
\accentedsymbol\tildepsi{\kern 3.2pt\tilde{\kern -3.2pt\psi}}
\accentedsymbol\ssizetildepsi{\ssize\kern 2.5pt\tilde{\kern -2.5pt\psi}}
%\accentedsymbol\hatgamma{\kern 2pt\hat{\kern -2pt\gamma}}
%\accentedsymbol\checkgamma{\kern 2.5pt\check{\kern -2.5pt\gamma}}
\def\Pover#1{\overset{\kern 1.5pt #1}\to P}
\def\Qover#1{\overset{\kern 1.5pt #1}\to Q}
\def\Rover#1{\overset{\kern 1.5pt #1}\to R}
\def\Sover#1{\overset{\kern 1.5pt #1}\to S}
\def\blue#1{#1}

\catcode`#=11\def\diez{#}\catcode`#=6
\catcode`&=11\catcode`&=4
\catcode`_=11\def\podcherkivanie{_}\catcode`_=8
\catcode`~=11\catcode`~=\active
\def\mycite#1{\cite{\blue{#1}}\immediate\special{ps:
     ShrHPSdict begin /ShrBORDERthickness 0 def}}
\def\myciterange#1#2#3#4{\cite{\blue{#2#3#4}}\immediate\special{ps:
     ShrHPSdict begin /ShrBORDERthickness 0 def}}
\def\mytag#1{%
    \tag#1}
\def\mythetag#1{\thetag{\blue{#1}}\immediate\special{ps:
     ShrHPSdict begin /ShrBORDERthickness 0 def}}
\def\myrefno#1{\no#1}
\def\myhref#1#2{\blue{#2}\immediate\special{ps:
     ShrHPSdict begin /ShrBORDERthickness 0 def}}
\def\myEarXivlink{\myhref{http://arXiv.org}{http:/\negskp/arXiv.org}}

\def\mytheorem#1{\csname proclaim\endcsname{Theorem #1}}
\def\mytheoremwithtitle#1#2{\csname proclaim\endcsname{Theorem #1#2}}
\def\mythetheorem#1{\blue{#1}\immediate\special{ps:
     ShrHPSdict begin /ShrBORDERthickness 0 def}}
\def\mylemma#1{\csname proclaim\endcsname{Lemma #1}}
\def\mylemmawithtitle#1#2{\csname proclaim\endcsname{Lemma #1#2}}
\def\mythelemma#1{\blue{#1}\immediate\special{ps:
     ShrHPSdict begin /ShrBORDERthickness 0 def}}
\def\mycorollary#1{\csname proclaim\endcsname{Corollary #1}}

\def\mydefinition#1{\definition{Definition #1}}
\def\mythedefinition#1{\blue{#1}\immediate\special{ps:
     ShrHPSdict begin /ShrBORDERthickness 0 def}}
\def\myconjecture#1{\csname proclaim\endcsname{Conjecture #1}}
\def\myconjecturewithtitle#1#2{\csname proclaim\endcsname{Conjecture #1#2}}

\def\myproblem#1{\csname proclaim\endcsname{Problem #1}}
\def\myproblemwithtitle#1#2{\csname proclaim\endcsname{Problem #1#2}}

%----------------------------------------------------------------
% Cyrillic fonts definition

%\font\tencyr=wncyr10
%----------------------------------------------------------------
\pagewidth{360pt}
\pageheight{606pt}
\topmatter
\title
Comparison of two classifications of a class of ODE's
in the first case of intermediate degeneration.
\endtitle
\rightheadtext{Comparison of two classifications.}
\author
Ruslan Sharipov
\endauthor
\address Bashkir State University, 32 Zaki Validi street, 450074 Ufa, Russia
\endaddress
\email\myhref{mailto:r-sharipov\@mail.ru}{r-sharipov\@mail.ru}
\endemail
\abstract
     Two classifications of second order ODE's cubic with respect to the 
first order derivative are compared in the first case of intermediate degeneration.
The correspondence of vectorial, pseudovectorial, scalar, and pseudoscalar invariants 
of the equations in this case is established. 
\endabstract
\endtopmatter
\loadbold
%\loadeufb
\TagsOnRight
\document
%\input countstyle

% \special{header=resource.eps}
\head
1. Introduction.
\endhead
     Since the epoch of classical papers (see \mycite{1} and \mycite{2}) it is
known that the class of second order differential equations cubic with respect 
to the first order derivative 
$$
\hskip -2em
y''=P(x,y)+3\,Q(x,y)\,y'+3\,R(x,y)\,(y')^2+S(x,y)\,(y')^3
\mytag{1.1}
$$
is closed with respect to transformations of the form 
$$
\hskip -2em
\cases
\tilde x=\tilde x(x,y),\\
\tilde y=\tilde y(x,y).
\endcases
\mytag{1.2}
$$
About 19 years ago in \mycite{3} and \mycite{4} the equations \mythetag{1.1} 
were classified using their scalar invariants. They were subdivided into 
{\bf nine subclasses} closed with respect to transformations of the form 
\mythetag{1.2}. The richest class comprising almost all equations of the form
\mythetag{1.1} consists of the equations of {\bf general position}. The smallest
class is composed by the equations of {\bf maximal degeneration}. The rest of
the equations \mythetag{1.1} are distributed among seven subclasses composed by 
the equations of {\bf intermediate degeneration}.\par
     Recently in 2013 Yu\.~Yu\.~Bagderina in \mycite{5} presented her own 
classification of the equations \mythetag{1.1} again subdividing them into 
{\bf nine subclasses} closed with respect to transformations of the form 
\mythetag{1.2}. She uses Sophus Lie's method of infinitesimal transformations 
adapted to equations of the form \mythetag{1.1} by N\.~H\.~Ibragimov 
in \mycite{6}.\par 
    In \mycite{5} Yu\.~Yu\.~Bagderina does not mention the previously existing 
classification from \mycite{3} and \mycite{4}. She cites the paper \mycite{4}
only as a source of invariants and for criticism of its method. The present 
paper is the second one in a series of papers intended to examine the results 
of \mycite{5} and compare them with the prior results from 
\myciterange{3}{3}{, }{4} and \mycite{7}. In the previous paper \mycite{8} it
was shown that items 1 and 9 in Bagderina's classification Theorem 2 in
\mycite{5} do coincide with the case of {\bf general position} and the case of
{\bf maximal degeneration} respectively from the previously existing 
classification in \mycite{3} and \mycite{4}. It was also revealed that in the 
case of {\bf general position} most structures and most formulas from
Bagderina's  paper \mycite{5} do coincide or are very closely related to those 
in \mycite{7}, though they are given in different notations.\par
     In the present paper we consider item 2 of Bagderina's classification 
Theorem in \mycite{5} and compare it with the {\bf first case of intermediate 
degeneration} in \mycite{3} and \mycite{4}. Then we study the 
structures and formulas from item 2 of Thheirem 2 in \mycite{5} and establish 
their correspondence to the structures and formulas of the previously existing 
classification in \mycite{3} and \mycite{4}.\par
\head
2. Some notations and definitions.
\endhead
    Transformations of the form \mythetag{1.2} are called point 
transformations. They are assumed to be locally invertible. The inverse 
transformations for them are also point transformations. They are written 
as follows:
$$
\hskip -2em
\cases
x=\tilde x(\tilde x,\tilde y),\\
y=\tilde y(\tilde x,\tilde y).
\endcases
\mytag{2.1}
$$
According to \myciterange{3}{3}{, }{4} and \mycite{7}, we use dot index
notations for partial derivatives, e\.\,g\. having two functions $f(x,y)$ and 
$g(\tilde x,\tilde y)$ we write
$$
\xalignat 2
\hskip -2em
f_{\ssize p.q}=\frac{\partial^{p+q}f}{\partial x^p\,\partial y^q},
&&g_{\ssize p.q}=\frac{\partial^{p+q}g}{\partial\tilde x^p\,
\partial\tilde y^q}.
\mytag{2.2}
\endxalignat
$$
In terms of the notations \mythetag{2.2} the Jacoby matrices of the direct
and inverse point transformations \mythetag{1.2} and \mythetag{2.1} are written 
as follows:
$$
\xalignat 2
&\hskip -2em
S=\Vmatrix
x_{\sssize 1.0} &x_{\sssize 0.1}\\
\vspace{1ex}
y_{\sssize 1.0} &y_{\sssize 0.1}
\endVmatrix,
&&T=\Vmatrix
\tilde x_{\sssize 1.0} &\tilde x_{\sssize 0.1}\\
\vspace{1ex}
\tilde y_{\sssize 1.0} &\tilde y_{\sssize 0.1}
\endVmatrix.
\mytag{2.3}
\endxalignat
$$
In geometry the transformations \mythetag{1.2} and \mythetag{2.1} are interpreted 
as changes of local curvilinear coordinates on the plane $\Bbb R^2$ or on some
two-dimensional manifold. Their Jacoby matrices are called the direct and inverse 
transition matrices (see \mycite{9}).\par
     Tensorial and pseudotensorial fields in local coordinates are presented
as arrays of functions whose arguments are $x,\,y$ or $\tilde x,\,\tilde y$ 
respectively. These arrays of functions are called their components. They obey
some definite transformation rules.
\mydefinition{2.1} A pseudotensorial field of the type $(r,s)$ and weight $m$ 
is an array of functions $F^{i_1\ldots\,i_r}_{j_1\ldots\,j_s}$ which under 
the change of coordinates \mythetag{1.2} transforms as 
$$
\hskip -2em
F^{i_1\ldots\,i_r}_{j_1\ldots\,j_s}=
(\det T)^m\sum\Sb p_1\ldots p_r\\ q_1\ldots q_s\endSb
S^{i_1}_{p_1}\ldots\,S^{i_r}_{p_r}\,\,
T^{q_1}_{j_1}\ldots\,T^{q_s}_{j_s}\,\,
\tilde F^{p_1\ldots\,p_r}_{q_1\ldots\,q_s}.
\mytag{2.4}
$$
Tensorial fields are those pseudotensorial fields whose weight $m$ in \mythetag{2.4}
is zero. The prefix ``pseudo'' always indicates the nonzero weight $m\neq 0$. 
\enddefinition
     Tensorial and pseudotensorial fields of the type $(1,0)$ are called vectorial 
and pseudovectorial fields. Tensorial and pseudotensorial fields of the type $(0,1)$ 
are called covectorial and pseudocovectorial fields. And finally, scalar and 
pseudoscalar fields are those fields whose type is $(0,0)$.
\mydefinition{2.2} Tensorial and pseudotensorial fields whose components are 
expressed through $y'$, through the coefficients $P$, $Q$, $R$, $S$ of the equation
\mythetag{1.1}, and through their partial derivatives are called tensorial and 
pseudotensorial invariants of this equation respectively.
\enddefinition
     {\bf Remark}. Typically, in differential geometry components of tensorial and
pseudotensorial fields depend on a point of the base manifold only, i\.\,e\. on $x$
and $y$ or on $\tilde x$ and $\tilde y$ in our particular case. If the dependence 
of other parameters is included, this makes an extension of the concept. So, in
Definition~\mythedefinition{2.2} we have extended tensorial and pseudotensorial 
fields.\par
\head
3. Some basic structures.
\endhead	
     In \mycite{5} Yu\.~Yu\. Bagderina introduces a long list of special 
notations. In order to distinguish her notations from those in 
\myciterange{3}{3}{, }{4} and \mycite{7} we use the upper 
mark {\tencyr<}Bgd{\tencyr>} for her notations. The first order expressions
introduced in \mycite{5} are 
$$
\hskip -2em
\aligned
&\alpha^{\sssize\text{Bgd}}_0=Q_{\sssize 1.0}-P_{\sssize 0.1}+2\,P\,R-2\,Q^2,\\
&\alpha^{\sssize\text{Bgd}}_1=R_{\sssize 1.0}-Q_{\sssize 0.1}+P\,S-Q\,R,\\
&\alpha^{\sssize\text{Bgd}}_2=S_{\sssize 1.0}-R_{\sssize 0.1}+2\,Q\,S-2\,R^2.
\endaligned
\mytag{3.1}
$$
The order of an expression is determined by the highest order of partial 
derivatives of $P$, $Q$, $R$, $S$ in it. As it was noted in \mycite{8},
Bagderina's alpha quantities \mythetag{3.1} coincide with the components 
of the symmetric two-dimensional array $\Omega$ from \mycite{7}:
$$
\xalignat 3
&\hskip -2em
\alpha^{\sssize\text{Bgd}}_0=\Omega_{11},
&&\alpha^{\sssize\text{Bgd}}_1=\Omega_{12}=\Omega_{21},
&&\alpha^{\sssize\text{Bgd}}_2=\Omega_{22}. 
\quad
\mytag{3.2}
\endxalignat
$$
The quantities $\Omega_{ij}$ in \mythetag{3.2} constitute nether a tensorial 
invariant nor a pseudotensor invariant. However, their derivatives are used
in constructing both tensorial and pseudotensor invariants.\par
     The second order expressions are given by the formulas \thetag{2.2} 
in \mycite{5}: 
$$
\hskip -2em
\aligned
&\beta^{\,\sssize\text{Bgd}}_1=\partial_x\alpha^{\sssize\text{Bgd}}_1
-\partial_y\alpha^{\sssize\text{Bgd}}_0+R\,\alpha^{\sssize\text{Bgd}}_0
-2\,Q\,\alpha^{\sssize\text{Bgd}}_1+P\,\alpha^{\sssize\text{Bgd}}_2,\\
&\beta^{\,\sssize\text{Bgd}}_2=\partial_x\alpha^{\sssize\text{Bgd}}_2
-\partial_y\alpha^{\sssize\text{Bgd}}_1+S\,\alpha^{\sssize\text{Bgd}}_0
-2\,R\,\alpha^{\sssize\text{Bgd}}_1+Q\,\alpha^{\sssize\text{Bgd}}_2.
\endaligned
\mytag{3.3}
$$
As it was noted in \mycite{8}, Bagderina's beta quantities \mythetag{3.3} 
coincide with the components of the pseudocovectorial field\/ $\boldsymbol\alpha$ 
of the weight $1$ constructed in \mycite{7}:
$$
\xalignat 2
&\hskip -2em
\beta^{\,\sssize\text{Bgd}}_1=\alpha_1=A,
&&\beta^{\,\sssize\text{Bgd}}_2=\alpha_2=B.
\mytag{3.4}
\endxalignat
$$
It is convenient to express the components \mythetag{3.4} of the field 
$\boldsymbol\alpha$ directly through $P$, $Q$, $R$, $S$, as it was done in 
\myciterange{3}{3}{, }{4} and \mycite{7}, rather than through \mythetag{3.1}:
$$
\hskip -2em
\aligned
&\aligned
 A=P_{\sssize 0.2}&-2\,Q_{\sssize 1.1}+R_{\sssize 2.0}+
 2\,P\,S_{\sssize 1.0}+S\,P_{\sssize 1.0}-\\
 \vspace{0.5ex}
 &-3\,P\,R_{\sssize 0.1}-3\,R\,P_{\sssize 0.1}-
 3\,Q\,R_{\sssize 1.0}+6\,Q\,Q_{\sssize 0.1},
 \endaligned\\
 \vspace{1ex}
&\aligned
 B=S_{\sssize 2.0}&-2\,R_{\sssize 1.1}+Q_{\sssize 0.2}-
 2\,S\,P_{\sssize 0.1}-P\,S_{\sssize 0.1}+\\
 \vspace{0.5ex}
 &+3\,S\,Q_{\sssize 1.0}+3\,Q\,S_{\sssize 1.0}+
 3\,R\,Q_{\sssize 0.1}-6\,R\,R_{\sssize 1.0}.
 \endaligned
\endaligned
\mytag{3.5}
$$\par
    Bagderina's third order expressions are given by the formulas 
\thetag{2.3} in \mycite{5}:  
$$
\hskip -2em
\aligned
&\gamma^{\sssize\text{Bgd}}_{10}=\partial_x\beta^{\,\sssize\text{Bgd}}_1
-Q\,\beta^{\,\sssize\text{Bgd}}_1+P\,\beta^{\,\sssize\text{Bgd}}_2,\\
&\gamma^{\sssize\text{Bgd}}_{11}=\partial_x\beta^{\,\sssize\text{Bgd}}_2
-R\,\beta^{\,\sssize\text{Bgd}}_1+Q\,\beta^{\,\sssize\text{Bgd}}_2,\\
&\gamma^{\sssize\text{Bgd}}_{20}=\partial_y\beta^{\,\sssize\text{Bgd}}_1
-R\,\beta^{\,\sssize\text{Bgd}}_1+Q\,\beta^{\,\sssize\text{Bgd}}_2,\\
&\gamma^{\sssize\text{Bgd}}_{21}=\partial_y\beta^{\,\sssize\text{Bgd}}_2
-S\,\beta^{\,\sssize\text{Bgd}}_1+R\,\beta^{\,\sssize\text{Bgd}}_2.
\endaligned
\mytag{3.6}
$$
The expressions \mythetag{3.6} are used in order to define other third
order expressions. They are given by the formulas \thetag{2.17} in
\mycite{5}:
$$
\hskip -2em
\aligned
&\Gamma^{\,\sssize\text{Bgd}}_0=3\,\beta^{\,\sssize\text{Bgd}}_2\,
\gamma^{\sssize\text{Bgd}}_{10}+\beta^{\,\sssize\text{Bgd}}_1\,
(\gamma^{\sssize\text{Bgd}}_{20}-4\,\gamma^{\sssize\text{Bgd}}_{11}),\\
&\Gamma^{\,\sssize\text{Bgd}}_1=\beta^{\,\sssize\text{Bgd}}_2
\,(4\,\gamma^{\sssize\text{Bgd}}_{20}-\gamma^{\sssize\text{Bgd}}_{11})
-3\,\beta^{\,\sssize\text{Bgd}}_1\,\gamma^{\sssize\text{Bgd}}_{21}.
\endaligned
\mytag{3.7}
$$
As it was shown in \mycite{8}, Bagderina's gamma quantities \mythetag{3.7} 
coincide with the components of the pseudocovectorial field\/ 
$\boldsymbol\beta$ of the weight $3$ constructed in \mycite{7}:
$$
\xalignat 2
&\hskip -2em
\Gamma^{\,\sssize\text{Bgd}}_0=\beta_1=-H,
&&\Gamma^{\,\sssize\text{Bgd}}_1=\beta_2=G.
\mytag{3.8}
\endxalignat
$$
The quantities $H$ and $G$ in \mythetag{3.8} can be expressed through 
$A$ and $B$ from \mythetag{3.5} in a more explicit way. They are given by 
the following formulas taken from \mycite{7}:
$$
\aligned
G&=-B\,B_{\sssize 1.0}-3\,A\,B_{\sssize 0.1}+4\,B\,
A_{\sssize 0.1}+3\,S\,A^2-6\,R\,B\,A+3\,Q\,B^2,\\
\vspace{1ex}
H&=-A\,A_{\sssize 0.1}-3\,B\,A_{\sssize 1.0}+4\,A\,
B_{\sssize 1.0}-3\,P\,B^2+6\,Q\,A\,B-3\,R\,A^2,
\endaligned
\quad
\mytag{3.9}
$$
\par
     Bagderina's fourth order expressions are given by the formulas 
\thetag{2.4} in \mycite{5}:  
$$
\hskip -2em
\aligned
&\delta^{\,\sssize\text{Bgd}}_{10}=\partial_x\gamma^{\sssize\text{Bgd}}_{10}
-2\,Q\,\gamma^{\sssize\text{Bgd}}_{10}+P\,(\gamma^{\sssize\text{Bgd}}_{20}
+\gamma^{\sssize\text{Bgd}}_{11})-5\,\alpha^{\sssize\text{Bgd}}_0\,
\beta^{\,\sssize\text{Bgd}}_1,\\
&\delta^{\,\sssize\text{Bgd}}_{20}=\partial_x\gamma^{\sssize\text{Bgd}}_{20}
-R\,\gamma^{\sssize\text{Bgd}}_{10}+P\,\gamma^{\sssize\text{Bgd}}_{21}
-4\,\alpha^{\sssize\text{Bgd}}_1\,\beta^{\,\sssize\text{Bgd}}_1
-\alpha^{\sssize\text{Bgd}}_0\,\beta^{\,\sssize\text{Bgd}}_2,\\
&\delta^{\,\sssize\text{Bgd}}_{30}=\partial_y\gamma^{\sssize\text{Bgd}}_{20}
-S\,\gamma^{\sssize\text{Bgd}}_{10}+Q\,\gamma^{\sssize\text{Bgd}}_{21}
-4\,\alpha^{\sssize\text{Bgd}}_2\,\beta^{\,\sssize\text{Bgd}}_1
-\alpha^{\sssize\text{Bgd}}_1\,\beta^{\,\sssize\text{Bgd}}_2,\\
&\delta^{\,\sssize\text{Bgd}}_{11}=\partial_x\gamma^{\sssize\text{Bgd}}_{11}
-R\,\gamma^{\sssize\text{Bgd}}_{10}+P\,\gamma^{\sssize\text{Bgd}}_{21}
-\alpha^{\sssize\text{Bgd}}_1\,\beta^{\,\sssize\text{Bgd}}_1
-4\,\alpha^{\sssize\text{Bgd}}_0\,\beta^{\,\sssize\text{Bgd}}_2,\\
&\delta^{\,\sssize\text{Bgd}}_{21}=\partial_x\gamma^{\sssize\text{Bgd}}_{21}
-R\,(\gamma^{\sssize\text{Bgd}}_{20}+\gamma^{\sssize\text{Bgd}}_{11})
+2\,Q\,\gamma^{\sssize\text{Bgd}}_{21}-5\,\alpha^{\sssize\text{Bgd}}_1
\,\beta^{\,\sssize\text{Bgd}}_2,\\
&\delta^{\,\sssize\text{Bgd}}_{31}=\partial_y\gamma^{\sssize\text{Bgd}}_{21}
-S\,(\gamma^{\sssize\text{Bgd}}_{20}+\gamma^{\sssize\text{Bgd}}_{11})
+2\,R\,\gamma^{\sssize\text{Bgd}}_{21}-5\,\alpha^{\sssize\text{Bgd}}_2
\,\beta^{\,\sssize\text{Bgd}}_2.
\endaligned
\mytag{3.10}
$$
The fifth order expressions by Bagderina are given by the formulas \thetag{2.5} 
in \mycite{5}:  
$$
\hskip -2em
\aligned
&\epsilon^{\,\sssize\text{Bgd}}_{10}=\partial_x\delta^{\,\sssize\text{Bgd}}_{10}
-3\,Q\,\delta^{\,\sssize\text{Bgd}}_{10}+P\,(2\,\delta^{\,\sssize\text{Bgd}}_{20}
+\delta^{\,\sssize\text{Bgd}}_{11})-12\,\alpha^{\sssize\text{Bgd}}_0\,
\gamma^{\sssize\text{Bgd}}_{10},\\
&\epsilon^{\,\sssize\text{Bgd}}_{20}=\partial_y\delta^{\,\sssize\text{Bgd}}_{10}
-3\,R\,\delta^{\,\sssize\text{Bgd}}_{10}+Q\,(2\,\delta^{\,\sssize\text{Bgd}}_{20}
+\delta^{\,\sssize\text{Bgd}}_{11})-12\,\alpha^{\sssize\text{Bgd}}_1\,
\gamma^{\sssize\text{Bgd}}_{10},\\
&\aligned
\epsilon^{\,\sssize\text{Bgd}}_{11}&=\partial_x\delta^{\,\sssize\text{Bgd}}_{11}
-R\,\delta^{\,\sssize\text{Bgd}}_{10}-Q\,\delta^{\,\sssize\text{Bgd}}_{11}
+2\,P\,\delta^{\,\sssize\text{Bgd}}_{21}-2\,\alpha^{\sssize\text{Bgd}}_1
\,\gamma^{\sssize\text{Bgd}}_{10}\,-\\
\vspace{-1ex}
&-\,10\,\alpha^{\sssize\text{Bgd}}_0\,\gamma^{\sssize\text{Bgd}}_{11}
-10\,(\beta^{\,\sssize\text{Bgd}}_1)^2.
\endaligned
\endaligned
\mytag{3.11}
$$
And finally, her sixth order expression is given by the formula \thetag{2.6} 
in \mycite{5}:  
$$
\hskip -2em
\lambda^{\,\sssize\text{Bgd}}_{10}=\partial_x\epsilon^{\,\sssize\text{Bgd}}_{10}
-4\,Q\,\epsilon^{\,\sssize\text{Bgd}}_{10}+P\,(3\,\epsilon^{\,\sssize\text{Bgd}}_{20}
+\epsilon^{\,\sssize\text{Bgd}}_{11})-21\,\alpha^{\sssize\text{Bgd}}_0
\,\delta^{\,\sssize\text{Bgd}}_{10}.
\mytag{3.12}
$$\par
     In $\Bbb R^2$ and in any two-dimensional manifold there are two pseudotensorial
fields with constant components. They are denoted by the same symbol $\bold d$ and
are given by the same skew-symmetric matrix in any local coordinates:
$$
\xalignat 2
&\hskip -2em
d_{ij}=\Vmatrix\format \r&\quad\l\\ 0 & 1\\-1 & 0\endVmatrix,
&&d^{\kern 1pt ij}=\Vmatrix\format \r&\quad\l\\ 0 & 1\\-1 & 0\endVmatrix.
\mytag{3.13}
\endxalignat
$$
The components $d_{ij}$ in \mythetag{3.13} correspond to the pseutotensorial
field $\bold d$ of the type $(0,2)$ and the weight $-1$. The components 
$d^{\kern 1pt ij}$ in \mythetag{3.13} correspond to the pseutotensorial 
field $\bold d$ of the type $(2,0)$ and the weight $1$. These two fields
are used for raising and lowering indices of other pseutotensorial fields. 
In particular, we have 
$$
\xalignat 2
&\hskip -2em
\alpha^i=\sum^2_{k=1}d^{\kern 1pt ik}\,\alpha_k,
&&\beta^i=\sum^2_{k=1}d^{\kern 1pt ik}\,\beta_k,
\mytag{3.14}
\endxalignat
$$ 
In explicit form the equalities \mythetag{3.14} are written as follows:
$$
\xalignat 2
&\hskip -2em
\alpha^1=B=\beta^{\,\sssize\text{Bgd}}_2,
&&\alpha^2=-A=-\beta^{\,\sssize\text{Bgd}}_1,
\mytag{3.15}\\
\vspace{1ex}
&\hskip -2em
\beta^1=G=\Gamma^{\,\sssize\text{Bgd}}_1,
&&\beta^2=H=-\Gamma^{\,\sssize\text{Bgd}}_0.
\mytag{3.16}
\endxalignat
$$
The quantities \mythetag{3.15} are the components of the pseudovectorial
field $\boldsymbol\alpha$ of the weight $2$. The quantities \mythetag{3.16} 
are the components of the pseudovectorial field $\boldsymbol\beta$ of the 
weight $4$. In \mycite{4} these quantities are used in order to define 
a pseudoscalar field $F$ of the weight $1$. This field is defined by means
of the formula
$$
\hskip -2em
3\,F^5=\sum^2_{i=1}\alpha_i\,\beta^i=-\sum^2_{i=1}\beta_i\,\alpha^i=
A\,G+B\,H.
\mytag{3.17}
$$
One can define $F$ more explicitly by applying \mythetag{3.5} and \mythetag{3.9}
to \mythetag{3.17}:
$$
\hskip -2em
\aligned
F^5=A\,B\,A_{\sssize 0.1}&+B\,A\,B_{\sssize 1.0}-
A^2\,B_{\sssize 0.1}-B^2\,A_{\sssize 1.0}-\\
&-P\,B^3+3\,Q\,A\,B^2-3\,R\,A^2\,B+S\,A^3.
\endaligned
\mytag{3.18}
$$
Yu\.~Yu\.~Bagderina introduces her own quantity $J_0$ (see \thetag{2.16} in
\mycite{5}):
$$
\hskip -2em
J^{\,\sssize\text{Bgd}}_0
=(\beta^{\,\sssize\text{Bgd}}_2)^2\,\gamma^{\sssize\text{Bgd}}_{10}
-\beta^{\,\sssize\text{Bgd}}_1\,\beta^{\,\sssize\text{Bgd}}_2\,
(\gamma^{\sssize\text{Bgd}}_{20}+\gamma^{\sssize\text{Bgd}}_{11})
+(\beta^{\,\sssize\text{Bgd}}_1)^2\,\gamma^{\sssize\text{Bgd}}_{21}.
\mytag{3.19}
$$
As it was noted in \mycite{8}, Bagderina's quantity $J^{\,\sssize\text{Bgd}}_0$
in \mythetag{3.19} is related to the quantity $F$ in \mythetag{3.17} and
in \mythetag{3.18} in the following way:
$$
\hskip -2em
J^{\,\sssize\text{Bgd}}_0=-F^5.
\mytag{3.20}
$$
Using $J^{\,\sssize\text{Bgd}}_0$, Yu\.~Yu\.~Bagderina introduces her
mu quantity $\mu^{\,\sssize\text{Bgd}}_1$:
$$
\pagebreak
\hskip -2em
\mu^{\,\sssize\text{Bgd}}_1=(J^{\,\sssize\text{Bgd}}_0)^{1/5}.
\mytag{3.21}
$$
Comparing the formula \mythetag{3.20} with the formula \mythetag{3.21}, we see 
that Bagderina's mu quantity differs from $F$ only in sign: 
$$
\mu^{\,\sssize\text{Bgd}}_1=-F.
$$
\par
     The pseudocovectorial and pseudovectorial fields $\boldsymbol\alpha$ 
with the components \mythetag{3.4} and \mythetag{3.15}, the pseudocovectorial 
and pseudovectorial fields $\boldsymbol\beta$ with the components \mythetag{3.8} 
and \mythetag{3.16}, and the pseudoscalar field $F$ in \mythetag{3.17} and
\mythetag{3.18} constitute basic structures associated with any equation of
the form \mythetag{1.1}. All of them are presented by Yu\.~Yu\.~Bagderina in 
\mycite{5} using her own notations. However, none of them is new in \mycite{5}
as compared to \myciterange{3}{3}{, }{4} and \mycite{7}.\par 
\head
4. Cases of intermediate degeneration.
\endhead
     According to the classification from \myciterange{3}{3}{, }{4} the class of
equations of the form \mythetag{1.1} is subdivided into {\bf nine subclasses} 
which are called cases. The case of {\bf general position} corresponds to the 
richest subclass of all nine. Any equation \mythetag{1.1} taken by chance falls 
into the case of general position with the probability 1. The other eight classes 
are thin classes. Their total measure (probability) is zero.\par
     The case of {\bf general position} is defined by the condition 
$$
\hskip -2em
F\neq 0. 
\mytag{4.1}
$$
Looking at \mythetag{3.17}, we can write \mythetag{4.1} as
$$
\hskip -2em
3\,F^5=\det\Vmatrix A & B\\\vspace{1.5ex} -H & G\endVmatrix\neq 0. 
\mytag{4.2}
$$
A matrix with nonzero determinant cannot have zero rows. It cannot have proportional
rows either. Note that $A$ and $B$ are components of the pseudocovectorial field
$\boldsymbol\alpha$ in \mythetag{3.4}, while $-H$ and $G$ are components of the pseudocovectorial field $\boldsymbol\beta$ in \mythetag{3.8}. Therefore, using
\mythetag{4.2}, the condition \mythetag{4.1} implies
$$
\xalignat 3
&\hskip -2em
\boldsymbol\alpha\neq 0,
&&\boldsymbol\beta\neq 0,
&&\boldsymbol\alpha\nparallel\boldsymbol\beta.
\mytag{4.3}
\endxalignat
$$
Conversely, using \mythetag{3.4}, \mythetag{3.8}, and \mythetag{4.2}, the 
conditions \mythetag{4.3} imply \mythetag{4.1}.\par
     The case of {\bf maximal degeneration} is opposite to \mythetag{4.3}.
It is given by the following condition for the pseudocovectorial field 
$\boldsymbol\alpha$:
$$
\hskip -2em
\boldsymbol\alpha=0.
\mytag{4.4}
$$
Using using \mythetag{3.4}, \mythetag{3.8}, and \mythetag{3.9}, the
condition \mythetag{4.4} implies $\boldsymbol\beta=0$ and $F=0$.\par
     Recently in \mycite{8} it was shown that the item 1 of Bagderina's
classification Theorem 2 in \mycite{5} is equivalent to the case of 
{\bf general position} from the prior papers \myciterange{3}{3}{, }{4} 
and \mycite{7}. Also in \mycite{8} it was shown that the item 9 
of this classification Theorem 2 in \mycite{5} is equivalent to the case of 
{\bf maximal degeneration} from \myciterange{3}{3}{, }{4}. Here we continue
our comparison work and proceed to the cases of {\bf intermediate 
degeneration}.\par
     The cases of {\bf intermediate degeneration} are splitted off from 
\mythetag{4.1} and \mythetag{4.4} by setting the following two conditions:
$$
\pagebreak
\xalignat 2
&\hskip -2em
F=0,
&&\boldsymbol\alpha\neq 0.
\mytag{4.5}
\endxalignat
$$
There are seven cases of {\bf intermediate degeneration}. In this paper we 
consider the first of them and compare it with item 2 in Bagderina's 
Theorem 2 in \mycite{5}.\par
    As we noted above, the conditions \mythetag{4.3} taken altogether
imply \mythetag{4.1}. Therefore, if the conditions \mythetag{4.5} are 
fulfilled, then at least one of the last two conditions in \mythetag{4.3}
should be broken. A zero vector (or a zero pseudovector) is parallel to
any other vector (or pseudovector). Therefore $\boldsymbol\beta=0$ is
just a subcase of a more general case $\boldsymbol\beta\parallel\boldsymbol
\alpha$. Hence the conditions \mythetag{4.5} are equivalent to the following
conditions:
$$
\xalignat 2
&\hskip -2em
\boldsymbol\alpha\neq 0,
&&\boldsymbol\beta\parallel\boldsymbol\alpha.
\mytag{4.6}
\endxalignat
$$
The conditions \mythetag{4.6} mean that there is a factor $N$ such that
$$
\hskip -2em
\boldsymbol\beta=3\,N\,\boldsymbol\alpha.
\mytag{4.7}
$$
The factor $N$ in \mythetag{4.7} is a pseudoscalar field of the weight $2$. 
The field $N$ was discovered in \mycite{3} in the form $N=Q$ in some special
coordinates (see \thetag{6.2} in \mycite{3}). The formula \mythetag{4.7} was 
presented in \mycite{4} (see \thetag{4.2} in \mycite{4}). The formulas 
$$
\xalignat 2
&\hskip -2em
N=\frac{G}{3\,B},&&N=-\frac{H}{3\,A}
\mytag{4.8}
\endxalignat
$$
are immediate from \mythetag{4.7}, \mythetag{4.8}, \mythetag{3.4}, and
\mythetag{3.8}, (see \thetag{4.3} in \mycite{4}). The first formula applies
in the case $B\neq 0$, the second one in the case $A\neq 0$. If both $A$ and
$B$ are nonzero, both formulas are applicable. Note that $A$ and $B$ cannot
vanish simultaneously since $A$ and $B$ are components of the field 
$\boldsymbol\alpha$ and $\boldsymbol\alpha\neq 0$ (see \mythetag{4.5} and
\mythetag{4.6} above).\par
{\bf Remark}. The pseudoscalar field $N$, i\.\,e\. an object with proper 
geometric behavior, arises only if the conditions \mythetag{4.5} are 
fulfilled. Otherwise the formulas \mythetag{4.8} yield two quantities with
no meaning at all.\par 
     Apart from $N$ there are some other pseudotensorial fields and objects 
of different geometric nature associated with the equations \mythetag{1.1} 
in all cases of intermediate degeneration. Here are the quantities
$\varphi_1$ and $\varphi_2$: 
$$
\hskip -2em
\aligned
&\varphi_1=-3\,A\,\frac{A\,S-B_{\sssize 0.1}}{5\,B^2}
-3\,\frac{A_{\sssize 0.1}+B_{\sssize 1.0}-3\,A\,R}{5\,B}
-\frac{6}{5}\,Q,\\
\vspace{2ex}
&\varphi_2=3\,\frac{A\,S-B_{\sssize 0.1}}{5\,B}
-\frac{3}{5}\,R.
\endaligned
\mytag{4.9}
$$
The formulas \mythetag{4.9} apply in the case $B\neq 0$. If $A\neq 0$, we
use the formulas 
$$
\hskip -2em
\aligned
&\varphi_1=
-3\,\frac{B\,P+A_{\sssize 1.0}}{5\,A}
+\frac{3}{5}\,Q,\\
\vspace{2ex}
&\varphi_2=
3\,B\,\frac{B\,P+A_{\sssize 1.0}}{5\,A^2}
-3\,\frac{B_{\sssize 1.0}+A_{\sssize 0.1}+3\,B\,Q}{5\,A}
+\frac{6}{5}\,R.
\endaligned
\mytag{4.10}
$$
If both $A$ and $B$ are nonzero, then both formulas \mythetag{4.9} and
\mythetag{4.10} are applicable.\par
     The quantities $\varphi_1$ and $\varphi_2$ do not form a pseudotensorial 
field. They are transformed as follows under the point transformations 
\mythetag{1.2}:
$$
\hskip -2em
\varphi_i=\sum^2_{j=1}T^j_i\,\tilde\varphi_j-
\frac{\partial\ln\det T}{\partial x^i}.
\mytag{4.11}
$$
Here $x^1=x$, $x^2=y$ and $T$ is the transition matrix defined in 
\mythetag{2.3}.\par
    The quantities $\varphi_1$ and $\varphi_2$ were introduced in \mycite{3}
using some special coordinates. The formulas \mythetag{4.9} and \mythetag{4.10}
were derived in \mycite{4}. These formulas are applicable in arbitrary 
coordinates $x$ and $y$. Due to the transformation rule \mythetag{4.11} the 
quantities $\varphi_1$ and $\varphi_2$ can be combined with the components 
of the array $\theta$ defined in \mycite{7} in order to form connection 
components (see \thetag{6.10} in \mycite{3} or \thetag{4.22} in \mycite{4}):
$$
\hskip -2em
\varGamma^k_{ij}=\theta^k_{ij}-\frac{\varphi_i\,\delta^k_j+
\varphi_j\,\delta^k_i}{3}.
\mytag{4.12}
$$
Note that the connection \mythetag{4.12} is different from the connection
used in \mycite{7} and later in \mycite{8} for the case of general 
position. The quantities $\varphi_1$ and $\varphi_2$ here are also 
different from those used in the case of general position.\par 
    The second field introduced in \mycite{3} is $M$ (see \thetag{6.7} in
\mycite{3}). The first was $N$ that was introduced as $N=Q$ in some special
coordinates (see \thetag{6.2} in \mycite{3}). The field $M$ was also
first introduced in that special coordinates. The formulas for $M$ in 
arbitrary coordinates were derived in \mycite{4} (see \thetag{4.28} and
\thetag{4.29} in \mycite{4}):
$$
\gather
\hskip -2em
\aligned
M=-\frac{12\,A\,N\,(A\,S-B_{\sssize 0.1})}{5\,B}
&-A\,N_{\sssize 0.1}+\frac{24}{5}\,A\,N\,R-\\
\vspace{1ex}
-\frac{6}{5}\,N\,A_{\sssize 0.1}-\frac{6}{5}\,N\,
&B_{\sssize 1.0}+B\,N_{\sssize 1.0}-\frac{12}{5}\,B\,N\,Q,
\endaligned
\mytag{4.13}\\
\vspace{2ex}
\hskip -2em
\aligned
M=-\frac{12\,B\,N\,(B\,P+A_{\sssize 1.0})}{5\,A}
&+B\,N_{\sssize 1.0}+\frac{24}{5}\,B\,N\,Q+\\
\vspace{1ex}
+\frac{6}{5}\,N\,B_{\sssize 1.0}+\frac{6}{5}\,N\,
&A_{\sssize 0.1}-A\,N_{\sssize 0.1}-\frac{12}{5}\,A\,N\,R.
\endaligned
\mytag{4.14}
\endgather
$$
The formula \mythetag{4.13} applies in the case $B\neq 0$. If $A\neq 0$, we
use the formula \mythetag{4.14}.\par
     In the cases of {\bf intermediate degeneration} we loose 
$\boldsymbol\beta$ as an independent field. It becomes parallel to 
$\boldsymbol\alpha$ (see \mythetag{4.6}). However, exactly at that instant
another pseudocovectorial field arises. It was discovered in \mycite{3} and
was denoted through $\boldsymbol\gamma$. Initially $\boldsymbol\gamma$ was
presented in some special coordinates. Then in \mycite{4} it was expressed
by explicit formulas in arbitrary coordinates (see \thetag{4.30} and 
\thetag{4.31} in \mycite{4}):
$$
\align
&\hskip -2em
\aligned
\gamma_1=&\frac{6\,A\,N\,(A\,S-B_{\sssize 0.1})}{5\,B^2}
-\frac{18\,N\,A\,R}{5\,B}+\\
\vspace{1ex}
&+\frac{6\,N\,(A_{\sssize 0.1}+B_{\sssize 1.0})}{5\,B}
-N_{\sssize 1.0}+\frac{12}{5}\,N\,Q-2\,\Omega\,A.
\endaligned
\quad
\mytag{4.15}\\
\vspace{2ex}
&\hskip -2em
\aligned
\gamma_2=-\frac{6\,N\,(A\,S-B_{\sssize 0.1})}{5\,B}
-N_{\sssize 0.1}+\frac{6}{5}\,N\,R-2\,\Omega\,B,
\endaligned
\quad
\mytag{4.16}
\endalign
$$
The formulas \mythetag{4.15} and \mythetag{4.16} are used if $B\neq 0$. 
If $A\neq 0$, we write:
$$
\align
&\hskip -2em
\aligned
\gamma_1=\frac{6\,N\,(B\,P+A_{\sssize 1.0})}{5\,A}
-N_{\sssize 1.0}-\frac{6}{5}\,N\,Q-2\,\Omega\,A,
\endaligned
\quad
\mytag{4.17}\\
\displaybreak
\vspace{2ex}
&\hskip -2em
\aligned
\gamma_2=-&\frac{6\,B\,N\,(B\,P+A_{\sssize 1.0})}{5\,A^2}
+\frac{18\,N\,B\,Q}{5\,A}+\\
\vspace{1ex}
&+\frac{6\,N\,(B_{\sssize 1.0}+A_{\sssize 0.1})}{5\,A}
-N_{\sssize 0.1}-\frac{12}{5}\,N\,R-2\,\Omega\,B.
\endaligned
\quad
\mytag{4.18}
\endalign
$$
The weight of the pseudocovectorial field $\boldsymbol\gamma$ given by 
the formulas \mythetag{4.15} and \mythetag{4.16} or by the formulas 
\mythetag{4.17} and \mythetag{4.18} is equal to $2$. Note that the
formulas \mythetag{4.15}, \mythetag{4.16}, \mythetag{4.17}, 
\mythetag{4.18} in \mycite{4} are given in a pseudovectorial form, 
i\.\,e\. with upper indices (see \thetag{4.30}, \thetag{4.31}, 
\thetag{4.32}, and \thetag{4.33} in \mycite{4}).\par
      The pseudoscalar field $\Omega$ is the third field common for all 
cases of {\bf intermediate degeneration}. This field was introduced in
\mycite{3} by means of the formulas
$$
\hskip -2em
\Omega=\frac{5}{6}\sum^2_{i=1}\sum^2_{j=1}\omega_{ij}\,d^{ij}
\text{, \ where \ }\omega_{ij}=\frac{\partial\varphi_i}{\partial x^j}
-\frac{\partial\varphi_j}{\partial x^i},
\mytag{4.19}
$$
in some special coordinates (see \thetag{6.17} and \thetag{6.18} in
\mycite{3}). It turns out that the formulas \mythetag{4.19} are applicable
in arbitrary coordinates as well (see \thetag{4.15} and \thetag{4.16} in
\mycite{4}). The matter is that $\Omega$ is related to the curvature tensor
of the connection \mythetag{4.12}. The well-known formula for the curvature 
tensor (see \mycite{9}) is written as
$$
\hskip -2em
R^k_{rij}=
\frac{\partial\varGamma^k_{jr}}{\partial x^i}
-\frac{\partial\varGamma^k_{ir}}{\partial x^j}+
\sum^2_{q=1}\varGamma^k_{iq}\varGamma^q_{jr}-
\sum^2_{q=1}\varGamma^k_{jq}\varGamma^q_{ir}.
\mytag{4.20}
$$
Like in \mythetag{4.11}, here $x^1=x$ and $x^2=y$ are coordinates. Substituting
\mythetag{4.12} into the formula \mythetag{4.19}, we calculate $R^k_{rij}$ and 
then find that 
$$
\hskip -2em
\omega_{ij}=\sum^2_{k=1}R^k_{kij}
\mytag{4.21}
$$
Due to \mythetag{4.21} the quantities $\omega_{ij}$ are components of a tensor,
while $\Omega$ in \mythetag{4.19} is a pseudoscalar of the weight $1$. Here
are explicit formulas for $\Omega$:
$$
\align
&\hskip -2em
\aligned
\Omega&=\frac{2\,A\,B_{\sssize 0.1}(A\,S-B_{\sssize 0.1})}{B^3}
+\frac{(2\,A_{\sssize 0.1}-3\,A\,R)\,B_{\sssize 0.1}}{B^2}+\\
\vspace{1ex}
&+\frac{(B_{\sssize 1.0}-2\,A_{\sssize 0.1})\,A\,S}{B^2}
+\frac{A\,B_{\sssize 0.2}-A^2\,S_{\sssize 0.1}}{B^2}-
\frac{A_{\sssize 0.2}}{B}+\\
\vspace{1ex}
&+\frac{3\,A_{\sssize 0.1}\,R+3\,A\,R_{\sssize 0.1}
-A_{\sssize 1.0}\,S-A\,S_{\sssize 1.0}}{B}
+R_{\sssize 1.0}-2\,Q_{\sssize 0.1},
\endaligned
\mytag{4.22}\\
\vspace{1ex}
&\hskip -2em
\aligned
\Omega&=\frac{2\,B\,A_{\sssize 1.0}(B\,P+A_{\sssize 1.0})}{A^3}
-\frac{(2\,B_{\sssize 1.0}+3\,B\,Q)\,A_{\sssize 1.0}}{A^2}+\\
\vspace{1ex}
&+\frac{(A_{\sssize 0.1}-2\,B_{\sssize 1.0})\,B\,P}{A^2}
-\frac{B\,A_{\sssize 2.0}+B^2\,P_{\sssize 1.0}}{A^2}+
\frac{B_{\sssize 2.0}}{A}+\\
\vspace{1ex}
&+\frac{3\,B_{\sssize 1.0}\,Q+3\,B\,Q_{\sssize 1.0}
-B_{\sssize 0.1}\,P-B\,P_{\sssize 0.1}}{A}
+Q_{\sssize 0.1}-2\,R_{\sssize 1.0}
\endaligned
\mytag{4.23}
\endalign
$$
(see \thetag{4.17} and \thetag{4.21} in \mycite{4}). The formula 
\mythetag{4.22} applies in the case $B\neq 0$. If $A\neq 0$, we
apply the formula \mythetag{4.23}.\par
     It is important to note that the fields $\boldsymbol\alpha$, 
$\boldsymbol\gamma$ and $M$ obey the relationship:
$$
\hskip -2em
M=\sum^2_{i=1}\alpha_i\,\gamma^i=\sum^2_{i=1}\sum^2_{j=1}\alpha_i\,
d^{\kern 1pt ij}\,\gamma_j.
\mytag{4.24}
$$
The relationship \mythetag{4.24} is easily derived from \thetag{4.27}
in \mycite{4}. Since $\boldsymbol\alpha$ and $\boldsymbol\gamma$ in the
right hand side of \mythetag{4.24} are pseudocovectorial fields of the
weights $1$ and $2$ respectively, and $\bold d$ is a pseudotensorial 
field of the weight $1$, we see that $M$ is a pseudoscalar field of the 
weight $4$. This fact is known since \mycite{3}.\par
\head
5. Special coordinates. 
\endhead
     Let's recall that the cases of {\bf intermediate degeneration}
were introduced and studied in \mycite{3} using some special coordinates
where $A=0$ and $B=1$. However, Bagderina's classification Theorem 2 in
\mycite{5} and her formulas are derived under the restriction 
$\beta^{\,\sssize\text{Bgd}}_1\neq 0$, which corresponds to $A\neq 0$
(see \mythetag{3.4} above). In order to compare our formulas with those
in Bagderina's paper \mycite{5} and in order to make this comparison
comfortable for us we need some other special coordinates, which are similar 
to those in \mycite{3}, but different from them. 
\mytheorem{5.1} For any equation \mythetag{1.1} with $\boldsymbol\alpha\neq 0$
in \mythetag{3.4} there are some variables $x$ and $y$ such that $A=1$ and 
$B=0$ in these variables. 
\endproclaim
\demo{Proof} Note that $\boldsymbol\alpha$ in \mythetag{3.4} is 
a peudocovectorial field of the weight $1$ associated with the equation 
\mythetag{1.1} through the formulas \mythetag{3.5}. Raising indices
according to \mythetag{3.14}, we get the pseudovectorial field 
$\boldsymbol\alpha$ of the weight $2$ with the components \mythetag{3.15}. 
Let $\bold X$ be some nonzero vector field such that 
$\bold X\parallel\boldsymbol\alpha$. We can choose such a field by fixing 
some coordinates $x$ and $y$ and setting $X^1=\alpha^1$ and $X^2=\alpha^2$ 
in these coordinates. Being fulfilled in some particular coordinates, due to
\mythetag{2.4} the parallelism $\bold X\parallel\boldsymbol\alpha$ holds 
in arbitrary coordinates.\par
     It is well-known that any vector field $\bold X$ can be straighten 
(see \mycite{10}), i\.\,e\. there are some coordinates $x$ and $y$ such that 
$$
\xalignat 2
&\hskip -2em
X^1=0,
&&X^2=1
\mytag{5.1}
\endxalignat
$$
in these coordinates. Combining \mythetag{5.1} with 
$\bold X\parallel\boldsymbol\alpha$ and $\boldsymbol\alpha\neq 0$, we get 
$$
\xalignat 2
&\hskip -2em
\alpha^1=B=0,
&&\alpha^2=-A\neq 0.
\mytag{5.2}
\endxalignat
$$
Now let's perform a special transformation of the form \mythetag{1.2} given by
the formulas
$$
\xalignat 2
&\hskip -2em
\tilde x=x
&&\tilde y=\tilde y(x,y).
\mytag{5.3}
\endxalignat
$$
For the transformed components of $\boldsymbol\alpha$ in \mythetag{5.2}
from \mythetag{2.4} and \mythetag{5.3} we derive 
$$
\hskip -2em
\Vmatrix 0\\ \vspace{2ex} -\tilde A\endVmatrix=
(\det T)^{-2}\,
\Vmatrix 1 & 0\\
\vspace{2ex}
\tilde y_{\sssize 1.0} & \tilde y_{\sssize 0.1}
\endVmatrix\cdot\Vmatrix 0\\ \vspace{2ex}
-A\endVmatrix.
\mytag{5.4}
$$
Applying \mythetag{2.3}, we find that \mythetag{5.4} is equivalent to
$\tilde B=0$ and $\tilde A=(\tilde y_{\sssize 0.1})^{-1}\,A$. It is clear 
that, choosing a proper function $\tilde y(x,y)$ in \mythetag{5.3}, we can 
reach the required equality $\tilde A=1$ in the transformed coordinates 
$\tilde x$ and $\tilde y$. 
\qed\enddemo
     Thus, due to Theorem~\mythetheorem{5.1} proved just above we have 
special coordinates such that the following equalities are fulfilled in them:
$$
\xalignat 2
&\hskip -2em
\alpha^1=B=0,
&&\alpha^2=-A=-1.
\mytag{5.5}
\endxalignat
$$
Assuming that such special coordinates are chosen for $x$ and $y$, we shall
apply \mythetag{5.5} to various formulas from previous sections.\par
     We cannot apply \mythetag{5.5} to the formulas \mythetag{4.8} since $B$ 
is in the denominator of the first of them. However we can apply \mythetag{5.5}
to \mythetag{4.7}. This yields
$$
\xalignat 2
&\hskip -2em
\beta^1=3\,N\,\alpha^1=0,
&&\beta^2=3\,N\,\alpha^2=-3\,N.
\mytag{5.6}
\endxalignat
$$  
Taking into account \mythetag{3.16}, from \mythetag{5.6} we derive 
$$
\xalignat 2
&\hskip -2em
G=\Gamma^{\,\sssize\text{Bgd}}_1=0,
&&H=-\Gamma^{\,\sssize\text{Bgd}}_0=-3\,N.
\mytag{5.7}
\endxalignat
$$     
On the other hand, substituting \mythetag{5.6} into the formulas
\mythetag{3.9}, we obtain
$$
\xalignat 2
&\hskip -2em
G=3\,S,
&&H=-3\,R.
\mytag{5.8}
\endxalignat
$$     
Comparing \mythetag{5.8} with \mythetag{5.7}, we find that in our special
coordinates
$$
\xalignat 2
&\hskip -2em
S=0,
&&N=R.
\mytag{5.9}
\endxalignat
$$\par
     The next step is to apply \mythetag{5.5} to \mythetag{4.23}. As a result
for the pseudoscalar field $\Omega$ we derive the following very simple 
formula:
$$
\hskip -2em
\Omega=Q_{\sssize 0.1}-2\,R_{\sssize 1.0}. 
\mytag{5.10}
$$
The pseudoscalar field $M$ is given by the formula \mythetag{4.14}. 
Applying \mythetag{5.5}, \mythetag{5.9}, and \mythetag{5.10} to this
formula, we obtain the following expression for $M$:
$$
\hskip -2em
M=-R_{\sssize 0.1}-\frac{12}{5}\,R^2.
\mytag{5.11}
$$
The components of the pseudocovectorial field $\boldsymbol\gamma$ are given by
the formulas \mythetag{4.17} and \mythetag{4.18}. Applying \mythetag{5.5},
\mythetag{5.9}, and \mythetag{5.10} to them and using \mythetag{5.11}, we get
$$
\xalignat 2
&\hskip -2em
\gamma_1=3\,R_{\sssize 1.0}-2\,Q_{\sssize 0.1}-\frac{6}{5}\,R\,Q,
&&\gamma_2=M.
\mytag{5.12}
\endxalignat
$$\par
     Unlike $\gamma_1$ and $\gamma_2$ in \mythetag{5.12}, the quantities 
$\varphi_1$ and $\varphi_2$ do not represent components of a pseudotensorial 
field. Nevertheless, applying \mythetag{5.5} to \mythetag{4.10}, we derive 
$$
\xalignat 2
&\hskip -2em
\varphi_1=\frac{3}{5}\,Q,
&&\varphi_2=\frac{6}{5}\,R.
\mytag{5.13}
\endxalignat
$$
The non-tensorial quantities \mythetag{5.13} are used in \mythetag{4.12} 
to define a connection.\par
     Now let's return to the section 3. The formulas \mythetag{3.1} and
\mythetag{3.2} for Bagderina's alpha quantities from \mycite{5} remain 
unchanged. The formulas \mythetag{3.4} express the following comparison 
lemma coinciding with Lemma~3.2 in \mycite{8}.
\mylemma{5.1} Bagderina's beta quantities \mythetag{3.3} coincide with the 
components of the pseudocovectorial field\/ $\boldsymbol\alpha$ of the weight 
$1$ constructed in \mycite{7}.
\endproclaim
\noindent
The formulas \mythetag{3.4} are affected by \mythetag{5.5} in our
special coordinates. They reduce to 
$$
\xalignat 2
&\hskip -2em
\beta^{\,\sssize\text{Bgd}}_1=1,
&&\beta^{\,\sssize\text{Bgd}}_2=0.
\mytag{5.14}
\endxalignat
$$     
The formulas \mythetag{3.3} are equivalent to \mythetag{3.5}. Due to 
\mythetag{5.14} or \mythetag{5.5} and due to $S=0$ in \mythetag{5.9} 
they lead to the following differential equations:
$$
\hskip -2em
\aligned
&\aligned
 P_{\sssize 0.2}&-2\,Q_{\sssize 1.1}+R_{\sssize 2.0}
 -3\,P\,R_{\sssize 0.1}\,-\\
 &-\,3\,R\,P_{\sssize 0.1}-
 3\,Q\,R_{\sssize 1.0}+6\,Q\,Q_{\sssize 0.1}=1,
 \endaligned\\
 \vspace{1.5ex}
&-2\,R_{\sssize 1.1}+Q_{\sssize 0.2}
 +3\,R\,Q_{\sssize 0.1}-6\,R\,R_{\sssize 1.0}=0.
\endaligned
\mytag{5.15}
$$
The equations \mythetag{5.15} can be used in order to express higher order 
derivatives through lower order ones.\par
     Bagderina's gamma quantities \mythetag{3.6} become very simple in our
special coordinates. They are given by the following formulas:
$$
\xalignat 2
&\hskip -2em
\gamma^{\sssize\text{Bgd}}_{10}=-Q,
&&\gamma^{\sssize\text{Bgd}}_{11}=-R,\\
\vspace{-1.5ex}
\mytag{5.16}\\
\vspace{-1.5ex}
&\hskip -2em
\gamma^{\sssize\text{Bgd}}_{20}=-R,
&&\gamma^{\sssize\text{Bgd}}_{21}=0.
\endxalignat
$$
Bagderina's gamma quantities \mythetag{3.7} are described by the following
comparison lemma coinciding with Lemma~3.5 in \mycite{8}.
\mylemma{5.2} Bagderina's gamma quantities \mythetag{3.7} coincide with the 
components of the pseudocovectorial field\/ $\boldsymbol\beta$ of the weight 
$3$ constructed in \mycite{7}.
\endproclaim
\noindent
Combining \mythetag{5.7} and \mythetag{5.8}, for Bagderina's gamma quantities 
\mythetag{3.7} in our special coordinates we derive the following formulas:
$$
\xalignat 2
&\hskip -2em
\Gamma^{\,\sssize\text{Bgd}}_0=3\,R,
&&\Gamma^{\,\sssize\text{Bgd}}_1=0.
\mytag{5.17}
\endxalignat
$$
The formulas \mythetag{5.17} are equally simple as \mythetag{5.16}. 
\par
     Apart from \mythetag{3.1}, \mythetag{3.3}, \mythetag{3.6},
\mythetag{3.10}, \mythetag{3.11}, and \mythetag{3.12}, 
Yu\.~Yu\.~Bagderina uses the quantity $J^{\,\sssize\text{Bgd}}_0$
in \mycite{5} (see \mythetag{3.19}). The quantity $J^{\,\sssize\text{Bgd}}_0$ 
in \mythetag{3.19} is described by the following comparison lemma coinciding 
with Lemma~3.3 in \mycite{8}.
\mylemma{5.3} Bagderina's quantity $J^{\,\sssize\text{Bgd}}_0$ from
\mythetag{3.19} is related to the pseudoscalar field $F$ of the weight\/ $1$
constructed in \mycite{7} by means of the formula
$$
\hskip -2em
J^{\,\sssize\text{Bgd}}_0=-F^5.
\mytag{5.18}
$$
\endproclaim 
In addition to \mythetag{3.19} Yu\.~Yu\.~Bagderina uses four other 
quantities in item 2 of her classification Theorem 2 in \mycite{5}. They 
are given by the formulas \thetag{2.17} in \mycite{5}:
$$
\pagebreak
j^{\,\sssize\text{Bgd}}_0=\frac{3}{\beta^{\,\sssize\text{Bgd}}_1}
\biggl(\frac{\beta^{\,\sssize\text{Bgd}}_2}{\beta^{\,\sssize\text{Bgd}}_1}
\,\delta^{\,\sssize\text{Bgd}}_{10}-\delta^{\,\sssize\text{Bgd}}_{11}
\biggr)+\frac{6\,\gamma^{\sssize\text{Bgd}}_{10}}
{(\beta^{\,\sssize\text{Bgd}}_1)^2}\biggl(\gamma^{\sssize\text{Bgd}}_{11}
-\frac{\beta^{\,\sssize\text{Bgd}}_2}{\beta^{\,\sssize\text{Bgd}}_1}
\,\gamma^{\sssize\text{Bgd}}_{11}\biggr).
\quad
\mytag{5.19}
$$
The quantity \mythetag{5.19} is the first of these additional quantities 
from \thetag{2.17} in \mycite{5}. The second one is given by the following 
formula:
$$
\aligned
j^{\,\sssize\text{Bgd}}_1&=\frac{5}{6}\biggl(2\,\beta^{\,\sssize\text{Bgd}}_2
\,\delta^{\,\sssize\text{Bgd}}_{20}-\beta^{\,\sssize\text{Bgd}}_1
\,\delta^{\,\sssize\text{Bgd}}_{30}-\frac{(\beta^{\,\sssize\text{Bgd}}_2)^2}
{\beta^{\,\sssize\text{Bgd}}_1}\,\delta^{\,\sssize\text{Bgd}}_{10}\biggr)
+\biggl(\gamma^{\sssize\text{Bgd}}_{20}\,-\\
&-\,\frac{2}{3}\,\gamma^{\sssize\text{Bgd}}_{11}
-\frac{\beta^{\,\sssize\text{Bgd}}_2}{3\,\beta^{\,\sssize\text{Bgd}}_1}
\,\gamma^{\sssize\text{Bgd}}_{10}\biggr)
\biggl(\gamma^{\sssize\text{Bgd}}_{20}+\gamma^{\sssize\text{Bgd}}_{11}
-2\,\frac{\beta^{\,\sssize\text{Bgd}}_2}{\beta^{\,\sssize\text{Bgd}}_1}
\,\gamma^{\sssize\text{Bgd}}_{10}\biggr).
\endaligned
\mytag{5.20}
$$
The rest two quantities are given by the formulas 
$$
\aligned
&\aligned
j^{\,\sssize\text{Bgd}}_2&=\frac{1}{\beta^{\,\sssize\text{Bgd}}_1}
\biggl(\delta^{\,\sssize\text{Bgd}}_{20}
-\frac{\beta^{\,\sssize\text{Bgd}}_2}{\beta^{\,\sssize\text{Bgd}}_1}
\,\delta^{\,\sssize\text{Bgd}}_{10}\biggr)\,+\\
&+\,\frac{\gamma^{\sssize\text{Bgd}}_{10}}{5\,(\beta^{\,\sssize\text{Bgd}}_1)^2}
\biggl(7\,\frac{\beta^{\,\sssize\text{Bgd}}_2}{\beta^{\,\sssize\text{Bgd}}_1}
\,\gamma^{\sssize\text{Bgd}}_{10}-6\,\gamma^{\sssize\text{Bgd}}_{20}
-\gamma^{\sssize\text{Bgd}}_{11}\biggr),
\endaligned\\
\vspace{1ex}
&j^{\,\sssize\text{Bgd}}_3=\frac{3}{5}\biggl(
\frac{\delta^{\,\sssize\text{Bgd}}_{10}}{(\beta^{\,\sssize\text{Bgd}}_1)^3}
-\frac{6\,(\gamma^{\sssize\text{Bgd}}_{10})^2}
{5\,(\beta^{\,\sssize\text{Bgd}}_1)^4}\biggr).
\endaligned
\mytag{5.21}
$$
The quantity $j^{\,\sssize\text{Bgd}}_0$ in \mythetag{5.19} is described by
the following comparison lemma. 
\mylemma{5.4} If the conditions \mythetag{4.5} are fulfilled, i\.\,e\. in the
cases of intermediate degeneration, Bagderina's jay quantity 
$j^{\,\sssize\text{Bgd}}_0$ from \mythetag{5.19} behaves as a pseudoscalar 
field of the weight $1$. It is related to the pseudoscalar field $\Omega$ 
introduced in \mycite{3} as
$$
\hskip -2em
j^{\,\sssize\text{Bgd}}_0=-3\,\Omega\,. 
\mytag{5.22}
$$
\endproclaim
Lemma~\mythelemma{5.4} is proved by verifying the formula \mythetag{5.22}.
This could be done directly using some symbolic algebra package. In my 
case that was Maple\footnotemark. 
\footnotetext{\ Maple is a trademark of Waterloo Maple Inc.}
\mylemma{5.5} If the conditions \mythetag{4.5} are fulfilled, i\.\,e\. in the
cases of intermediate degeneration, Bagderina's jay quantity 
$j^{\,\sssize\text{Bgd}}_1$ from \mythetag{5.20} behaves as a pseudoscalar 
field of the weight $4$. It is related to the pseudoscalar field $M$ 
introduced in \mycite{3} as
$$
\hskip -2em
j^{\,\sssize\text{Bgd}}_1=\frac{5}{2}\,M. 
\mytag{5.23}
$$
\endproclaim
    Lemma~\mythelemma{5.5} is similar to Lemma~\mythelemma{5.4}. It is proved
by verifying the formula \mythetag{5.23} by means of direct computations.\par
    Bagderina's quantities $j^{\,\sssize\text{Bgd}}_2$ and 
$j^{\,\sssize\text{Bgd}}_3$ in \mythetag{5.21} are different. They do not 
behave as pseudoscalar fields. However, some definite combination of them do. 
On page 27 of her paper \mycite{5} Yu\.~Yu\.~Bagderina introduces the following 
quantity:
$$
\hskip -2em
j^{\,\sssize\text{Bgd}}_5=5\,\bigl(2\,j^{\,\sssize\text{Bgd}}_1\,
j^{\,\sssize\text{Bgd}}_3+(j^{\,\sssize\text{Bgd}}_2
-j^{\,\sssize\text{Bgd}}_0/6)^2\bigr).
\mytag{5.24}
$$
\mylemma{5.6} If the conditions \mythetag{4.5} are fulfilled, i\.\,e\. in the
cases of intermediate degeneration, Bagderina's jay quantity 
$j^{\,\sssize\text{Bgd}}_5$ from \mythetag{5.24} behaves as a pseudoscalar 
\pagebreak field of the weight $2$. 
\endproclaim
     The relation of Bagderina's field $j^{\,\sssize\text{Bgd}}_5$ to the fields 
introduced in \mycite{3} and \mycite{4} is studied below. Now we write the explicit 
formulas for $j^{\,\sssize\text{Bgd}}_2$ and $j^{\,\sssize\text{Bgd}}_3$ from
\mythetag{5.21} in our special coordinates introduced according to
Theorem~\mythetheorem{5.1}:
$$
\hskip -2em
\aligned
&j^{\,\sssize\text{Bgd}}_2=4\,Q_{\sssize 0.1}-5\,R_{\sssize 1.0}
+\frac{18}{5}\,Q\,R,\\
\vspace{1ex}
&j^{\,\sssize\text{Bgd}}_3=3\,P_{\sssize 0.1}-\frac{18}{5}\,Q_{\sssize 1.0}
-\frac{36}{5}\,P\,R+\frac{162}{25}\,Q^2.
\endaligned
\mytag{5.25}
$$
The formula for $j^{\,\sssize\text{Bgd}}_5$ in these special coordinates
is more complicated than \mythetag{5.25}:
$$
\hskip -2em
\gathered
j^{\,\sssize\text{Bgd}}_5=180\,R\,P\,R_{\sssize 0.1}
-216\,Q\,R\,R_{\sssize 1.0}
-\bigl(180\,R_{\sssize 0.0}^2+75\,R_{\sssize 0.1}\bigr)\,P_{\sssize 0.1}\,+\\
\vspace{1ex}
+\,\bigl(216\,R^2+90\,R_{\sssize 0.1}\bigr)\,Q_{\sssize 1.0}
-\bigl(270\,R_{\sssize 1.0}
-162\,Q\,R\bigr)\,Q_{\sssize 0.1}\,-\\
\vspace{1ex}
-\,162\,R_{\sssize 0.1}\,Q^2
+180\,R_{\sssize 1.0}^2
+432\,R^3\,P-324\,Q^2\,R^2+\frac{405}{4}\,Q_{\sssize 0.1}^2.
\endgathered
\mytag{5.26}
$$
\par
     Apart from $j^{\,\sssize\text{Bgd}}_5$, on page 27 of her paper
\mycite{5} Yu\.~Yu\.~Bagderina introduces the quantity 
$j^{\,\sssize\text{Bgd}}_4$ by means of the following formula:
$$
\hskip -2em
j^{\,\sssize\text{Bgd}}_4=\frac{\Gamma^{\,\sssize\text{Bgd}}_0}
{\beta^{\,\sssize\text{Bgd}}_1}.
\mytag{5.27}
$$
\mylemma{5.7} If the conditions \mythetag{4.5} are fulfilled, i\.\,e\. in the
cases of intermediate degeneration, under the auxiliary condition 
$\beta^{\,\sssize\text{Bgd}}_1\neq 0$ Bagderina's jay quantity 
$j^{\,\sssize\text{Bgd}}_4$ from \mythetag{5.27} is related to the 
pseudoscalar field $N$ of the weight $2$ introduced in \mycite{3} and 
effectively calculated in \mycite{4} by means of the formula
$$
\hskip -2em
j^{\,\sssize\text{Bgd}}_4=3\,N. 
\mytag{5.28}
$$ 
\endproclaim
     The comparison Lemma~\mythelemma{5.7} is immediate from 
Lemma~\mythelemma{5.1} and Lemma~\mythelemma{5.2} due to the formulas 
\mythetag{3.4}, \mythetag{3.8}, and \mythetag{4.8}. 
\head
6. The first case of intermediate degeneration
and Bagderina's type two equations. 
\endhead
     The {\bf first case of intermediate degeneration} is determined
by the conditions
$$
\xalignat 3
&\hskip -2em
F=0,
&&\boldsymbol\alpha\neq 0,
&&M\neq 0,
\mytag{6.1}
\endxalignat
$$
where $F=0$ and $\boldsymbol\alpha\neq 0$ are common for all cases
{\bf of intermediate degeneration} (see \mythetag{4.5}). From $M\neq 0$,
using either \mythetag{4.13} or \mythetag{4.14}, we derive 
$$
\hskip -2em
N\neq 0,
\mytag{6.2}
$$
where $N$ is given by the formulas \mythetag{4.8}. The equality 
\mythetag{4.24} can be written as
$$
\hskip -2em
M=\det\Vmatrix \alpha_1 & \alpha_2\\
\vspace{1.5ex}\gamma_1 & \gamma_2\endVmatrix, 
\mytag{6.3}
$$
which is similar to \mythetag{4.2}. From $M\neq 0$ and \mythetag{6.3} 
we derive 
$$
\xalignat 3
&\hskip -2em
\boldsymbol\alpha\neq 0,
&&\boldsymbol\gamma\neq 0,
&&\boldsymbol\alpha\nparallel\boldsymbol\gamma.
\mytag{6.4}
\endxalignat
$$
Conversely, due to \mythetag{6.3} the conditions \mythetag{6.4} imply 
$M\neq 0$, i\.\,e\. they are equivalent to the inequality $M\neq 0$. 
\par
    Bagderina's type two equations are defined in item 2 of her
classification Theorem~2 in \mycite{5}. They are given by the
following conditions:
$$
\xalignat 4
&\hskip -2em
J^{\,\sssize\text{Bgd}}_0\neq 0,
&&\beta^{\,\sssize\text{Bgd}}_1\neq 0,
&&j^{\,\sssize\text{Bgd}}_0\neq 0,
&&\Gamma^{\,\sssize\text{Bgd}}_0\neq 0.
\quad
\mytag{6.5}
\endxalignat
$$
Applying the comparison lemmas (see Lemma~\mythelemma{5.3}, 
Lemma~\mythelemma{5.1}, Lemma~\mythelemma{5.4}, Lem\-ma~\mythelemma{5.7}
and the formulas \mythetag{5.18}, \mythetag{3.15}, \mythetag{5.22},
\mythetag{5.27}, \mythetag{5.28}), we can write the conditions
\mythetag{6.5} in terms of the fields introduced in \myciterange{3}{3}{, }{4} 
and \mycite{7}:
$$
\xalignat 4
&\hskip -2em
F\neq 0,
&&\boldsymbol\alpha\neq 0,
&&\Omega\neq 0,
&&N\neq 0.
\quad
\mytag{6.6}
\endxalignat
$$
Comparing \mythetag{6.6} with \mythetag{6.1}, we see that the conditions
do not coincide. This means that Bagderina's classification in \mycite{5} 
is slightly different from that of \myciterange{3}{3}{, }{4}. For the further 
comparison purposes we draw the following table.\par
\medskip
\newdimen\contentsize
\contentsize=\hsize
\advance\contentsize by -3pt
\newdimen\halfcontentsize
\halfcontentsize=\contentsize
\advance\halfcontentsize by -7pt
\divide\halfcontentsize by 2
\newdimen\onethirdcontentsize
\onethirdcontentsize=\contentsize
\advance\onethirdcontentsize by -10pt
\divide\onethirdcontentsize by 3
\newdimen\onethirdcontentsizeplus
\onethirdcontentsizeplus=\onethirdcontentsize
\advance\onethirdcontentsizeplus by 1pt
\newdimen\onethirdcontentsizeminus
\onethirdcontentsizeminus=\onethirdcontentsize
\advance\onethirdcontentsizeminus by 0pt
\newdimen\onethirdcontentsizeplusplus
\onethirdcontentsizeplus=\onethirdcontentsize
\advance\onethirdcontentsizeplusplus by -4.5pt
\newdimen\questionsize
\newdimen\biletwidth
\def\vtrule{\vrule height 13pt depth 5pt}
%----------------------------------------
\def\mytableline#1#2{%
\vtrule
\hbox to\contentsize{
\hbox to\onethirdcontentsize{\hss #1\hss}
\vtrule
\hbox to\onethirdcontentsize{\hss \hss}
\vtrule
\hbox to\onethirdcontentsize{\hss #2\hss}
\hss}\hskip 2pt plus 6pt minus 0pt
\vtrule\newline
\hbox to\contentsize{
\vbox{\hsize=\onethirdcontentsize\advance\hsize by 4.5pt
\hrule width\hsize}
\kern\onethirdcontentsizeminus
%\vbox{\hsize=\onethirdcontentsizeminus
%\hbox to \hsize{}}
\vbox{\hsize=\onethirdcontentsize\advance\hsize by 2pt
\hrule width\hsize}\hss}\newline
}
%----------------------------------------
\vbox{\baselineskip 0pt \offinterlineskip
\hrule width\hsize
\noindent
\vtrule
\hbox to\contentsize{
\hbox to\halfcontentsize{\hss R\.~A\.~Sharipov's classification\hss}
\vtrule
\hbox to\halfcontentsize{\hss Yu\.~Yu\.~Bagderina's classification\hss}
\hss}\hskip 2pt plus 3pt minus 0pt
\vtrule
\newline
\vrule height 8pt depth 5pt
\hbox to\contentsize{
\hbox to\halfcontentsize{\hss 1997-1998\hss}
\vrule height 8pt depth 5pt
\hbox to\halfcontentsize{\hss 2013\hss}
\hss}\hskip 2pt plus 3pt minus 0pt
\vrule height 8pt depth 5pt
\hrule width\hsize
\noindent
\vtrule
\hbox to\contentsize{
\hbox to\halfcontentsize{\hss ShrGP\hss}
\vtrule
\hbox to\halfcontentsize{\hss BgdET1\hss}
\hss}\hskip 2pt plus 3pt minus 0pt
\vtrule
\hrule width\hsize
\noindent
\mytableline{ShrID1}{BgdET2}%
\mytableline{ShrID2}{BgdET3}%
\mytableline{ShrID3}{BgdET4}%
\mytableline{ShrID4}{BgdET5}%
\mytableline{ShrID5}{BgdET6}%
\mytableline{ShrID6}{BgdET7}%
\vtrule
\hbox to\contentsize{
\hbox to\onethirdcontentsize{\hss ShrID7\hss}
\vtrule
\hbox to\onethirdcontentsize{\hss  \hss}
\vtrule
\hbox to\onethirdcontentsize{\hss BgdET8\hss}
\hss}\hskip 2pt plus 6pt minus 0pt
\vtrule\newline
\hrule width\hsize
\noindent
\vtrule
\hbox to\contentsize{
\hbox to\halfcontentsize{\hss ShrMD\hss}
\vtrule
\hbox to\halfcontentsize{\hss BgET9\hss}
\hss}\hskip 2pt plus 3pt minus 0pt
\vtrule
\hrule width\hsize
\vskip 0.2cm
}
\medskip
\noindent
The abbreviations in the above table read as follows:
\roster
\rosteritemwd=5pt
\item"--" ShrGP stands for Sharipov's case of general position;
\item"--" ShrMD stands for Sharipov's case of maximal degeneration;
\item"--" ShrID1 stands for Sharipov's case of intermediate degeneration 1;
\item"  " \. \. \. \. \. \. \. \. \. \. \. \. \. \. \. \. \. \. \. \. \. \.
          \. \. \. \. \. \. \. \. \. \. \. \. \. \. \. \. \. \. \. \. \. \.
          \. \. \.
\item"--" ShrID1 stands for Sharipov's case of intermediate degeneration 7; 
\item"--" BgdET1 stands for Bagderina's equations of type 1;
\item"  " \. \. \. \. \. \. \. \. \. \. \. \. \. \. \. \. \. \. \. \. \. \.
          \. \. \. \. \. \. \. \. \. \. \. \. \. \. \.
\item"--" BgdET9 stands for Bagderina's equations of type 9.
\endroster
\par
     Again looking at \mythetag{6.6} and \mythetag{6.1}, we see that generally
speaking the equation classes ShrID1 and BgdET2 do not coincide, but they have
a substantial overlap. Their overlap is described by the following conditions:
$$
\xalignat 4
&\hskip -2em
F\neq 0,
&&\boldsymbol\alpha\neq 0,
&&M\neq 0,
&&\Omega\neq 0.
\quad
\mytag{6.7}
\endxalignat
$$
Indeed, $M\neq 0$ in \mythetag{6.1} implies $N\neq 0$ in \mythetag{6.6}
(see \mythetag{6.2}). However $M\neq 0$ in \mythetag{6.1} does not imply 
$\Omega\neq 0$ in \mythetag{6.6}, unless some deeper mutual relations of 
these field will be discovered. Conversely, $N\neq 0$ and $\Omega\neq 0$ 
in \mythetag{6.6} do not imply $M\neq 0$ in \mythetag{6.1}. Below we shall 
study the intersection class ShrID1\,$\cap$\,BgdET2.\par
     Let's consider Bagderina's invariant differentiation operator
$\Cal D^{\,\sssize\text{Bgd}}_1$. It is given by the first formula 
\thetag{2.8} from Bagderina's classification Theorem 2 in \mycite{5}: 
$$
\hskip -2em
\Cal D^{\,\sssize\text{Bgd}}_1=\frac{\beta^{\,\sssize\text{Bgd}}_2}
{(\mu^{\,\sssize\text{Bgd}}_1)^2}
\,\frac{\partial}{\partial x}
-\frac{\beta^{\,\sssize\text{Bgd}}_1}
{(\mu^{\,\sssize\text{Bgd}}_1)^2}\,\frac{\partial}{\partial y}.
\mytag{6.8}
$$
Taking into account \mythetag{3.15}, \mythetag{5.22} and Bagderina's formula
$\mu^{\,\sssize\text{Bgd}}_1=j^{\,\sssize\text{Bgd}}_0$ from \thetag{2.10}
in item 2 of Theorem 2 in \mycite{5} (which is different from \mythetag{3.21}), 
we write \mythetag{6.8} as 
$$
\hskip -2em
\Cal D^{\,\sssize\text{Bgd}}_1=\frac{\alpha^1}{(3\,\Omega)^2}
\,\frac{\partial}{\partial x}
+\frac{\alpha^2}{(3\,\Omega)^2}\,\frac{\partial}{\partial y}.
\mytag{6.9}
$$
Invariant differentiation operators were not considered in \mycite{3} and 
\mycite{4} for the first case of intermediate degeneration ShrID1. Instead
of them covariant differentiation operators along pseudovectorial fields
were considered. In particular we have
$$
\xalignat 2
&\hskip -2em
\nabla_{\boldsymbol\alpha}=\alpha^1\,\nabla_1+\alpha^2\,\nabla_2,
&&\nabla_{\boldsymbol\gamma}=\gamma^1\,\nabla_1+\gamma^2\,\nabla_2
\mytag{6.10}
\endxalignat
$$
(see \thetag{6.13} in \mycite{3} and \thetag{5.2} in \mycite{4}). The covariant 
derivatives in \mythetag{6.10} extend partial derivatives from \mythetag{6.9}.
They are defined by means of the formula
$$
\hskip -2em
\aligned
\nabla_kF^{i_1\ldots\,i_r}_{j_1\ldots\,j_s}&=
\frac{\partial F^{i_1\ldots\,i_r}_{j_1\ldots\,j_s}}{\partial x^k}+
\sum^r_{n=1}\sum^2_{v_n=1}\varGamma^{i_n}_{k\,v_n}\,
F^{i_1\ldots\,v_n\ldots\,i_r}_{j_1\ldots\,j_s}-\\
&-\sum^s_{n=1}\sum^2_{w_n=1}\varGamma^{w_n}_{k\,j_n}\,
F^{i_1\ldots\,i_r}_{j_1\ldots\,w_n\ldots\,j_s}+
m\,\varphi_k\,F^{i_1\ldots\,i_r}_{j_1\ldots\,j_s}
\endaligned
\mytag{6.11}
$$
(see \thetag{6.11} in \mycite{3} or \thetag{4.23} in \mycite{4}). The 
covariant derivative $\nabla_k$ in \mythetag{6.11} is applied to a 
pseudotensorial field of the type $(r,s)$ and the weight $m$. The connection 
components $\varGamma^k_{ij}$ in \mythetag{6.11} are defined by \mythetag{4.12}. 
They are canonically associated with a given equation equation \mythetag{1.1}.
\par
     Due to the covariant derivatives $\nabla_1$ and $\nabla_2$ in \mythetag{6.10}
the differential operators \mythetag{6.10} are applicable not only to scalar invariants,
but to any tensorial and pseudotensorial invariants as well. Yu\.~Yu\.~Bagderina's
operator \mythetag{6.9} can also be extended in this manner using covariant 
derivatives $\nabla_1$ and $\nabla_2$:
$$
\hskip -2em
\Cal D^{\,\sssize\text{Bgd}}_1=\frac{\alpha^1}{(3\,\Omega)^2}
\,\nabla_1+\frac{\alpha^2}{(3\,\Omega)^2}\,\nabla_2.
\mytag{6.12}
$$
\mylemma{6.1} Within the intersection class\/ {\rm ShrID1\,$\cap$\,BgdET2}, 
i\.\,e\. if the conditions \mythetag{6.7} are fulfilled, Bagderina's 
invariant differentiation operator $\Cal D^{\,\sssize\text{Bgd}}_1$ from
\mythetag{6.8} extended in \mythetag{6.12} is related to the 
differentiation operator $\nabla_{\boldsymbol\alpha}$ from \mycite{3} as 
$$
\hskip -2em
\Cal D^{\,\sssize\text{Bgd}}_1=\frac{1}{{(3\,\Omega)^2}}
\,\nabla_{\boldsymbol\alpha}.
\mytag{6.13}
$$
\endproclaim
     Lemma~\mythelemma{6.1} and the formula \mythetag{6.13} are immediate from 
\mythetag{6.9} and \mythetag{6.10}. So we can proceed to the second invariant 
differentiation operator by Yu\.~Yu\.~Bagderina. It is given by the second 
formula \thetag{2.8} in Bagderina's Theorem 2 in \mycite{5}: 
$$
\hskip -2em
\Cal D^{\,\sssize\text{Bgd}}_2=
\biggl(\mu^{\,\sssize\text{Bgd}}_2\,\beta^{\,\sssize\text{Bgd}}_2
-3\,\frac{\mu^{\,\sssize\text{Bgd}}_2}{\beta^{\,\sssize\text{Bgd}}_1}
\biggr)\,\frac{\partial}{\partial x}
+\beta^{\,\sssize\text{Bgd}}_1\,\frac{\partial}{\partial y}.
\mytag{6.14}
$$
The quantity $\mu^{\,\sssize\text{Bgd}}_2$ in \mythetag{6.14} is given
by one of the formulas \thetag{2.10} in item 2 of Bagderina's classification 
Theorem 2 in \mycite{5}: 
$$
\hskip -2em
\mu^{\,\sssize\text{Bgd}}_2=\frac{3\,\beta^{\,\sssize\text{Bgd}}_1
\,e^{\,\sssize\text{Bgd}}_1}{\Gamma^{\,\sssize\text{Bgd}}_0}.
\mytag{6.15}
$$
The quantity $e^{\,\sssize\text{Bgd}}_1$ in \mythetag{6.15} is expressed
by one of the formulas \thetag{2.18} in \mycite{5}: 
$$
\hskip -2em
\gathered
e^{\,\sssize\text{Bgd}}_1=\frac{5}{(\beta^{\,\sssize\text{Bgd}}_1)^2}
\biggl(\frac{\beta^{\,\sssize\text{Bgd}}_2}{\beta^{\,\sssize\text{Bgd}}_1}
\,\epsilon^{\,\sssize\text{Bgd}}_{10}-\epsilon^{\,\sssize\text{Bgd}}_{11}
\biggr)\,+\\
\vspace{1ex}
+\,\frac{15}{(\beta^{\,\sssize\text{Bgd}}_1)^3}\biggl(
\gamma^{\sssize\text{Bgd}}_{11}-\frac{\beta^{\,\sssize\text{Bgd}}_2}
{\beta^{\,\sssize\text{Bgd}}_1}\,\gamma^{\sssize\text{Bgd}}_{10}\biggr)
-\frac{6\,\gamma^{\sssize\text{Bgd}}_{10}}{(\beta^{\,\sssize\text{Bgd}}_1)^2}
\,j^{\,\sssize\text{Bgd}}_0.
\endgathered
\mytag{6.16}
$$
The notations used in \mythetag{6.16} are given above in \mythetag{3.3},
\mythetag{3.6}, \mythetag{3.10}, \mythetag{3.11} and in \mythetag{5.19}. 
In order to transform \mythetag{6.14} we use our special coordinates 
introduced according to Theorem~\mythetheorem{5.1}. Upon replacing partial
derivatives in \mythetag{6.14} by covariant derivatives $\nabla_1$ and 
$\nabla_2$ in these special coordinates we get
$$
\aligned
\Cal D^{\,\sssize\text{Bgd}}_2&=
\bigl(9\,Q_{\sssize 0.1}-18\,R_{\sssize 1.0}\bigr)\,\nabla_1+\biggl(
\frac{10\,P_{\sssize 0.2}}{R}-30\,P_{\sssize 0.1}\,-\\
\vspace{1ex}
&-\,\frac{15\,Q_{\sssize 1.1}}{R}
-\frac{36\,Q\,R_{\sssize 1.0}}{R}+\frac{63\,Q\,Q_{\sssize 0.1}}{R}
-\frac{30\,P\,R_{\sssize 0.1}}{R}
-\frac{60}{R}\biggr)\,\nabla_2.
\endaligned
\quad
\mytag{6.17}
$$
In order to reveal the invariant nature of the operator \mythetag{6.17}
we calculate the covariant derivatives of the pseudoscalar field $\Omega$
in our special coordinates:
$$
\aligned
&\hskip -2em
\aligned
\nabla_1\Omega=
2\,P_{\sssize 0.2}
&-3\,Q_{\sssize 1.1}-6\,R\,P_{\sssize 0.1}\,-\\
\vspace{1ex}
&-\,\frac{36\,Q\,R_{\sssize 1.0}}{5}+\frac{63\,Q\,Q_{\sssize 0.1}}{5}
-6\,P\,R_{\sssize 0.1}-2,
\endaligned\\
\vspace{2ex}
&\hskip -2em
\nabla_2\Omega=\frac{18\,R\,R_{\sssize 1.0}}{5}-\frac{9\,R\,Q_{\sssize 0.1}}{5}.
\endaligned
\mytag{6.18}
$$
The quantities \mythetag{6.18} are components of the pseudocovectorial field
$\nabla\Omega$ of the weight 1. In order to apply them to \mythetag{6.17} we
need to raise their indices:
$$
\hskip -2em
\nabla^i\Omega=\sum^2_{k=1}d^{\kern 1pt ik}\,\nabla_k\Omega\,.
\mytag{6.19}
$$
The quantities \mythetag{6.19} are components of the pseudovectorial field
$\nabla\Omega$ of the weight 2. Due to \mythetag{3.13} the formula \mythetag{6.19}
simplifies to 
$$
\xalignat 2
&\hskip -2em
\nabla^1\Omega=\nabla_2\Omega,
&&\nabla^2\Omega=-\nabla_1\Omega\,.
\mytag{6.20}
\endxalignat
$$
Taking into account \mythetag{5.5} and \mythetag{6.20}, then comparing \mythetag{6.18} 
with \mythetag{6.17}, we get
$$
\hskip -2em
\Cal D^{\,\sssize\text{Bgd}}_2=\biggl(\frac{50\,\alpha^1}{N}
-\frac{5\,\nabla^1\Omega}{N}\biggr)\,\nabla_1
+\biggl(\frac{50\,\alpha^2}{N}
-\frac{5\,\nabla^2\Omega}{N}\biggr)\,\nabla_2.
\mytag{6.21}
$$
The denominator $R$ in \mythetag{6.17} is replaced by the denominator
$N$ in \mythetag{6.21} since $N=R$ in our special coordinates (see
\mythetag{5.9}).\par
     Note that \mythetag{6.21} is a proper tensorial formula. Therefore,
being derived in our special coordinates, it remains valid in arbitrary 
coordinates.\par
     Any vectorial and/or pseudovectorial field on the plane $\Bbb R^2$
or in a two-dimen\-sional manifold can be expressed as a linear combination
of any other two non-parallel vectorial and/or pseudovectorial field. In
our case ShrID1\,$\cap$\,BgdET2 this means that $\nabla\Omega$
is expressed through $\boldsymbol\alpha$ and $\boldsymbol\gamma$ since
$\boldsymbol\alpha\nparallel\boldsymbol\gamma$ due to $M\neq 0$ in
\mythetag{6.7}. This expression for $\nabla\Omega$ can be written 
explicitly:
$$
\hskip -2em
\nabla\Omega=\frac{\nabla_{\boldsymbol\gamma}\Omega}{M}\,\boldsymbol\alpha
-\frac{\nabla_{\boldsymbol\alpha}\Omega}{M}\,\boldsymbol\gamma.
\mytag{6.22}
$$
Now, applying \mythetag{6.22} to \mythetag{6.21}, we derive the following
formula:
$$
\hskip -2em
\Cal D^{\,\sssize\text{Bgd}}_2=
\biggl(\frac{50}{N}-
\frac{5\,\nabla_{\boldsymbol\gamma}\Omega}{M\,N}\biggr)
\nabla_{\boldsymbol\alpha}
+\biggl(\frac{5\,\nabla_{\boldsymbol\alpha}\Omega}{M\,N}\biggr)
\,\nabla_{\boldsymbol\gamma}.
\mytag{6.23}
$$
\mylemma{6.2} Within the intersection class\/ {\rm ShrID1\,$\cap$\,BgdET2}, 
i\.\,e\. if the conditions \mythetag{6.7} are fulfilled, Bagderina's 
invariant differentiation operator $\Cal D^{\,\sssize\text{Bgd}}_2$ in 
\mythetag{6.14} is related to the covariant differentiation operators 
$\nabla_{\boldsymbol\alpha}$ and $\nabla_{\boldsymbol\gamma}$ introduced in
\mycite{3} according to the formula \mythetag{6.23}.
\endproclaim     
\head
7. Curvature tensor and additional fields.
\endhead
    Let's return to the field $\Omega$ associated with the curvature
tensor \mythetag{4.20} of the connection \mythetag{4.20}. Following
the receipt of \mycite{3} we write
$$
\hskip -2em
R^k_{qij}=R^k_q\,d_{ij}\text{, \ where \ }
R^k_q=\frac{1}{2}\sum^2_{i=1}\sum^2_{j=1}R^k_{qij}\,d^{ij}
\mytag{7.1}
$$
(see \thetag{7.1} and \thetag{7.2} in \mycite{3}). The quantities $R^k_q$ 
in \mythetag{7.1} are components of a pseudotensorial field of the type 
$(1,1)$ and the the weight $1$. This field has two pseudoscalar invariants
--- its trace and its determinant. \pagebreak The trace of this field 
reduces to the pseudoscalar field $\Omega$ according to the formula: 
$$
\hskip -2em
\tr(R)=\frac{3}{5}\,\Omega. 
\mytag{7.2}
$$
Its determinant is a new field. In the framework of the {\bf second case of
intermediate degenerations} (see \mycite{4}), i\.\,e\. if $M=0$, this field 
can be expressed through the field $\Lambda$ given by the formulas 
\thetag{6.10} in \mycite{4}:
$$\
\hskip -2em
\det(R)=-\frac{9}{25}\,\Lambda\,(\Omega+\Lambda).
\mytag{7.3}
$$
In the present paper we deal with the case $M\neq 0$. Therefore we shall not
use the formulas \mythetag{7.2} and \mythetag{7.3} and we shall treat $\det(R)$ 
as an separate pseudoscalar field of the weight $2$. This field can be easily 
calculated in arbitrary coordinates using the formulas \mythetag{4.12}, 
\mythetag{4.20}, and \mythetag{7.1}. However, we choose our special coordinates 
introduced through Theorem~\mythetheorem{5.1}. In these coordinates we have 
$$
\gathered
\det(R)=
-\frac{36}{35}\,R\,P\,R_{\sssize 0.1}
+\frac{216}{125}\,Q\,R\,R_{\sssize 1.0}
+\biggl(\frac{36}{25}\,R_{\sssize 0.0}^2+\frac{3}{5}\,R_{\sssize 0.1}\biggr)
\,P_{\sssize 0.1}\,-\\
\vspace{1ex}
-\,\biggl(\frac{216}{125}\,R^2
+\frac{18}{25}\,R_{\sssize 0.1}\biggr)\,Q_{\sssize 1.0}
+\biggl(\frac{9}{5}\,R_{\sssize 1.0}
-\frac{162}{125}\,Q\,R\biggr)\,Q_{\sssize 0.1}\,+\\
\vspace{1ex}
+\,\frac{162}{125}\,R_{\sssize 0.1}\,Q^2
-\frac{27}{25}\,R_{\sssize 1.0}^2
-\frac{432}{125}\,R^3\,P
+\frac{324}{125}\,Q^2\,R^2
-\frac{18}{25}\,Q_{\sssize 0.1}^2.
\endgathered
\mytag{7.4}
$$
Comparing \mythetag{7.4} with \mythetag{5.26}, we derive the following formula:
$$
\hskip -2em
j^{\,\sssize\text{Bgd}}_5=-125\,\det(R)+\frac{45}{4}\,\Omega^2.
\mytag{7.5}
$$
This formula \mythetag{7.5} proves Lemma~\mythelemma{5.6}. It expresses 
Bagderina's quantity $j^{\,\sssize\text{Bgd}}_5$ through the pseudotensorial
field $R$ in \mythetag{7.1} previously known in \mycite{4}.\par
     Our further efforts are toward expressing $\det(R)$ through $M$, $N$,
$\Omega$, $\boldsymbol\alpha$, $\boldsymbol\gamma$ and their proper tensorial
derivatives. For this purpose we need a little bit of theory.\par
     For a while assume that $\boldsymbol\alpha$ and $\boldsymbol\gamma$ are 
arbitrary two pseudovectorial fields with the weights $m$ and $n$ respectively. 
Let $\bold X$ be a third pseudovectorial field with the weight $k$. Then we 
have the following identities:
$$
\align
&\hskip -2em
[\nabla_{\boldsymbol\alpha},\nabla_{\boldsymbol\gamma}]\bold X
-\nabla_{[\boldsymbol\alpha,\boldsymbol\gamma]}\bold X
=\bold R(\boldsymbol\alpha,\boldsymbol\gamma)\bold X
-k\ \omega(\boldsymbol\alpha,\boldsymbol\gamma)\ \bold X,
\mytag{7.6}\\
\vspace{1ex}
&\hskip -2em
\nabla_{\boldsymbol\alpha}\boldsymbol\gamma-
\nabla_{\boldsymbol\gamma}\boldsymbol\alpha
=[\boldsymbol\alpha,\boldsymbol\gamma]+
\bold T(\boldsymbol\alpha,\boldsymbol\gamma)\bold X.
\mytag{7.7}
\endalign
$$
Here $\bold R(\boldsymbol\alpha,\boldsymbol\gamma)$ and
$\bold T(\boldsymbol\alpha,\boldsymbol\gamma)$ are the curvature operator
and the torsion operator respectively (see \mycite{11}). The term 
$k\ \omega(\boldsymbol\alpha,\boldsymbol\gamma)\ \bold X$ is determined
by the skew symmetric form $\omega$ whose components are given in
\mythetag{4.19}. The connection components \mythetag{4.12} are symmetric.
Therefore we have no torsion here:
$$
\hskip -2em
\bold T(\boldsymbol\alpha,\boldsymbol\gamma)=0.
\mytag{7.8}
$$\par
     The formulas \mythetag{7.6} and \mythetag{7.7} are well known in
differential geometry, though they are usually applied to vectorial fields
rather then to pseudovectorial ones. Their application to pseudovectorial
fields have some features. In particular, the covariant derivatives 
\mythetag{6.11} and the commutator of pseudovectorial fields requires 
some auxiliary quantities $\varphi_i$ obeying the transformation rules
\mythetag{4.11}:
$$
\hskip -2em
[\boldsymbol\alpha,\boldsymbol\gamma]=\sum^2_{i=1}
\biggl(\sum^2_{s=1}\alpha^s\,\frac{\partial\gamma^i}{\partial x^s}
-\gamma^s\,\frac{\partial\alpha^i}{\partial x^s}
+n\,\alpha^s\,\varphi_s\,\gamma^i-m\,\gamma^s\,\varphi_s\,\alpha^i
\biggl)\frac{\partial}{\partial x^i}.
\mytag{7.9}
$$
The formula \mythetag{7.9} can be treated as a definition of commutator
in the case of pseudotensorial fields.\par
    Returning to our previously defined pseudovectorial fields 
$\boldsymbol\alpha$ and $\boldsymbol\gamma$ we should remind that 
their weights are $2$ and $3$ respectively. Note that they were 
originally defined as pseudocovectorial fields of the weights $1$ 
and $2$. But having the skew-symmetric metric pseudotensors 
\mythetag{3.13}, we can always raise and lower indices of any 
pseudotensorial field. We should also note that 
$$
\hskip -2em
\nabla\bold d=0
\mytag{7.10}
$$
for both metric pseudotensors with the components \mythetag{3.13}. Due
to \mythetag{7.10} the operations of raising and lowering indices 
commute with covariant differentiations.\par
     Now let's calculate the curvature operator $\bold R(\boldsymbol\alpha,
\boldsymbol\gamma)$ applied to some pseudovectorial field $\bold X$ taking
into account the special structure of $R^k_{qij}$ in \mythetag{7.1}:
$$
\hskip -2em
\bold R(\boldsymbol\alpha,\boldsymbol\gamma)\bold X
=\sum^2_{i=1}\sum^2_{j=1}\sum^2_{s=1}\sum^2_{q=1}
R^s_q\,d_{ij}\,\alpha^i\,\gamma^j\,X^q\,\frac{\partial}{\partial x^s}.
\mytag{7.11}
$$
Here $R^s_q$ are the components of that very matrix whose determinant is
applied in \mythetag{7.5}. Taking into account \mythetag{4.24}, we write
\mythetag{7.11} as 
$$
\hskip -2em
\bold R(\boldsymbol\alpha,\boldsymbol\gamma)\bold X
=\sum^2_{s=1}\sum^2_{q=1}
M\,R^s_q\,X^q\,\frac{\partial}{\partial x^s}=M\,R(\bold X).
\mytag{7.12}
$$
Here $R$ is the linear operator whose matrix is formed by $R^s_q$.\par
     At this moment we can apply \mythetag{7.6} to \mythetag{7.12}. As
a result we get 
$$
\hskip -2em
R(\bold X)=\frac{1}{M}\Bigl([\nabla_{\boldsymbol\alpha},\nabla_{\boldsymbol\gamma}]
\bold X-\nabla_{[\boldsymbol\alpha,\boldsymbol\gamma]}\bold X
+k\ \omega(\boldsymbol\alpha,\boldsymbol\gamma)\ \bold X\Bigr).
\mytag{7.13}
$$
Let's recall the formula \mythetag{4.21}. This formula combined
with \mythetag{7.1} and \mythetag{4.24} yields
$$
\hskip -2em
\omega(\boldsymbol\alpha,\boldsymbol\gamma)=\tr(R)\,
\sum^2_{i=1}\sum^2_{j=1}d_{ij}\,\alpha^i\,\gamma^j=\tr(R)\,M.
\mytag{7.14}
$$
Substituting \mythetag{7.14} into \mythetag{7.13}, we derive 
$$
\hskip -2em
R(\bold X)=\frac{1}{M}\Bigl([\nabla_{\boldsymbol\alpha},\nabla_{\boldsymbol\gamma}]
\bold X-\nabla_{[\boldsymbol\alpha,\boldsymbol\gamma]}\bold X\Bigr)
+k\,\tr(R)\ \bold X.
\mytag{7.15}
$$
The commutator $[\boldsymbol\alpha,\boldsymbol\gamma]$ in \mythetag{7.15} can
be calculated with the use of \mythetag{7.7} and  \mythetag{7.8}: 
$$
\hskip -2em
[\boldsymbol\alpha,\boldsymbol\gamma]=\nabla_{\boldsymbol\alpha}\boldsymbol\gamma
-\nabla_{\boldsymbol\gamma}\boldsymbol\alpha. 
\mytag{7.16}
$$
The covariant derivatives $\nabla_{\boldsymbol\alpha}\boldsymbol\gamma$ and
$\nabla_{\boldsymbol\gamma}\boldsymbol\alpha$ are that very derivatives used
in \myciterange{3}{3}{, }{4}: 
$$
\xalignat 2
&\hskip -2em
\nabla_{\boldsymbol\alpha}\boldsymbol\alpha=\Gamma^1_{11}\,\boldsymbol\alpha+
\Gamma^2_{11}\,\boldsymbol\gamma,
&&\nabla_{\boldsymbol\alpha}\boldsymbol\gamma=\Gamma^1_{12}\,\boldsymbol\alpha+
\Gamma^2_{12}\,\boldsymbol\gamma,\\
\vspace{-1.7ex}
&&&\mytag{7.17}\\
\vspace{-1.7ex}
&\hskip -2em
\nabla_{\boldsymbol\gamma}\boldsymbol\alpha=\Gamma^1_{21}\,\boldsymbol\alpha+
\Gamma^2_{21}\,\boldsymbol\gamma,
&&\nabla_{\boldsymbol\gamma}\boldsymbol\gamma=\Gamma^1_{22}\,\boldsymbol\alpha+
\Gamma^2_{22}\,\boldsymbol\gamma
\endxalignat
$$
(see \thetag{6.13} in \mycite{3} or \thetag{5.2} in \mycite{4}). 
The coefficients $\Gamma^1_{11}$, $\Gamma^2_{11}$, $\Gamma^1_{12}$, 
$\Gamma^2_{12}$, $\Gamma^1_{21}$, $\Gamma^2_{21}$, $\Gamma^1_{22}$, 
$\Gamma^2_{22}$ are pseudoscalar fields uniquely determined by 
$\boldsymbol\alpha$ and $\boldsymbol\gamma$ since $\boldsymbol\alpha
\nparallel\boldsymbol\gamma$ due to $M\neq 0$ in \mythetag{6.7}. 
Applying \mythetag{7.17} to \mythetag{7.16}, we get
$$
\hskip -2em
[\boldsymbol\alpha,\boldsymbol\gamma]
=(\Gamma^1_{12}-\Gamma^1_{21})\,\boldsymbol\alpha
+(\Gamma^2_{12}-\Gamma^2_{21})\,\boldsymbol\gamma
\mytag{7.18}
$$
Then we use \mythetag{7.18} in order to calculate $\nabla_{[\boldsymbol\alpha,
\boldsymbol\gamma]}\bold X$ in \mythetag{7.15}:
$$
\hskip -2em
\nabla_{[\boldsymbol\alpha,\boldsymbol\gamma]}\bold X
=(\Gamma^1_{12}-\Gamma^1_{21})\,\nabla_{\boldsymbol\alpha}\bold X
+(\Gamma^2_{12}-\Gamma^2_{21})\,\nabla_{\boldsymbol\gamma}\bold X. 
\mytag{7.19}
$$
The first term $[\nabla_{\boldsymbol\alpha}\nabla_{\boldsymbol\gamma}]
\bold X$ in \mythetag{7.15} is expanded as follows:
$$
\hskip -2em
[\nabla_{\boldsymbol\alpha},\nabla_{\boldsymbol\gamma}]
\bold X=\nabla_{\boldsymbol\alpha}\bigl(\nabla_{\boldsymbol\gamma}
\bold X\bigr)-\nabla_{\boldsymbol\gamma}\bigl(\nabla_{\boldsymbol\alpha}
\bold X\bigr).
\mytag{7.20}
$$
Taking into account \mythetag{7.19} and \mythetag{7.20}, the formula 
\mythetag{7.15} turns to
$$
\hskip -2em
\gathered
R(\bold X)=\frac{\nabla_{\boldsymbol\alpha}\bigl(\nabla_{\boldsymbol\gamma}
\bold X\bigr)-\nabla_{\boldsymbol\gamma}\bigl(\nabla_{\boldsymbol\alpha}
\bold X\bigr)}{M}\,-\\
\vspace{2ex}
-\frac{(\Gamma^1_{12}-\Gamma^1_{21})}{M}\,\nabla_{\boldsymbol\alpha}\bold X
-\frac{(\Gamma^2_{12}-\Gamma^2_{21})}{M}\,\nabla_{\boldsymbol\gamma}\bold X
+k\,\tr(R)\ \bold X.
\endgathered
\mytag{7.21}
$$
\par
     Now, using the formula \mythetag{7.21}, we can calculate two pseudovectoial
fields $R(\boldsymbol\alpha)$ and $R(\boldsymbol\gamma)$, expressing them back 
through $\boldsymbol\alpha$ and $\boldsymbol\gamma$ as linear combinations:
$$
\xalignat 2
&\hskip -2em
R(\boldsymbol\alpha)=F^1_1\,\boldsymbol\alpha+F^2_1\,\boldsymbol\gamma,
&&R(\boldsymbol\gamma)=F^1_2\,\boldsymbol\alpha+F^2_2\,\boldsymbol\gamma.
\mytag{7.22}
\endxalignat
$$
Indeed, for $R(\boldsymbol\alpha)$ in \mythetag{7.22} we have the following 
expression (with $k=2$):
$$
\hskip -2em
\gathered
R(\boldsymbol\alpha)=\frac{\nabla_{\boldsymbol\alpha}\bigl(\nabla_{\boldsymbol\gamma}
\boldsymbol\alpha\bigr)-\nabla_{\boldsymbol\gamma}\bigl(\nabla_{\boldsymbol\alpha}
\boldsymbol\alpha\bigr)}{M}\,-\\
\vspace{2ex}
-\frac{(\Gamma^1_{12}-\Gamma^1_{21})}{M}\,\nabla_{\boldsymbol\alpha}\boldsymbol\alpha
-\frac{(\Gamma^2_{12}-\Gamma^2_{21})}{M}\,\nabla_{\boldsymbol\gamma}\boldsymbol\alpha
+k\,\tr(R)\ \boldsymbol\alpha.
\endgathered
\mytag{7.23}
$$
The derivatives $\nabla_{\boldsymbol\alpha}\boldsymbol\alpha$ and 
$\nabla_{\boldsymbol\gamma}\boldsymbol\alpha$ are taken from \mythetag{7.17}. 
The second order derivatives $\nabla_{\boldsymbol\alpha}\bigl(\nabla_{\boldsymbol
\gamma}\boldsymbol\alpha\bigr)$ and $\nabla_{\boldsymbol\alpha}
\bigl(\nabla_{\boldsymbol\gamma}\boldsymbol\alpha\bigr)$ are transformed as follows:
$$
\gather
\hskip -2em
\gathered
\nabla_{\boldsymbol\alpha}\bigl(\nabla_{\boldsymbol
\gamma}\boldsymbol\alpha\bigr)=
\nabla_{\boldsymbol\alpha}\bigl(\Gamma^1_{21}\,\boldsymbol\alpha+
\Gamma^2_{21}\,\boldsymbol\gamma\bigr)=
\nabla_{\boldsymbol\alpha}\bigl(\Gamma^1_{21}\bigr)\,\boldsymbol\alpha\,+\\
\vspace{0.5ex}
+\,\Gamma^1_{21}\,\bigl(\Gamma^1_{11}\,\boldsymbol\alpha+
\Gamma^2_{11}\,\boldsymbol\gamma\bigr)+\nabla_{\boldsymbol\alpha}\bigl(
\Gamma^2_{21}\bigr)\,\boldsymbol\gamma+\Gamma^2_{21}
\,\bigl(\Gamma^1_{12}\,\boldsymbol\alpha+\Gamma^2_{12}\,\boldsymbol\gamma\bigr),
\endgathered
\quad
\mytag{7.24}\\
\vspace{1ex}
\hskip -2em
\gathered
\nabla_{\boldsymbol\gamma}\bigl(\nabla_{\boldsymbol
\alpha}\boldsymbol\alpha\bigr)=\nabla_{\boldsymbol\gamma}\bigl(
\Gamma^1_{11}\,\boldsymbol\alpha+\Gamma^2_{11}\,\boldsymbol\gamma\bigr)
=\nabla_{\boldsymbol\gamma}\bigl(\Gamma^1_{11}\bigr)\,\,\boldsymbol\alpha\,+\\
\vspace{0.5ex}
+\,\Gamma^1_{11}\,\bigl(\Gamma^1_{21}\,\boldsymbol\alpha+
\Gamma^2_{21}\,\boldsymbol\gamma\bigr)
+\nabla_{\boldsymbol\gamma}\bigl(\Gamma^2_{11}\bigr)\,\boldsymbol\gamma
+\Gamma^2_{11}\,\bigl(\Gamma^1_{22}\,\boldsymbol\alpha
+\Gamma^2_{22}\,\boldsymbol\gamma\bigr).
\endgathered
\quad
\mytag{7.25}
\endgather
$$
The pseudovectorial field $R(\boldsymbol\gamma)$ in \mythetag{7.22} is 
treated similarly. From \mythetag{7.21} we derive the following expression 
for this field (with $k=3$):
$$
\hskip -2em
\gathered
R(\boldsymbol\gamma)=\frac{\nabla_{\boldsymbol\alpha}\bigl(\nabla_{\boldsymbol\gamma}
\boldsymbol\gamma\bigr)-\nabla_{\boldsymbol\gamma}\bigl(\nabla_{\boldsymbol\alpha}
\boldsymbol\gamma\bigr)}{M}\,-\\
\vspace{2ex}
-\frac{(\Gamma^1_{12}-\Gamma^1_{21})}{M}\,\nabla_{\boldsymbol\alpha}\boldsymbol\gamma
-\frac{(\Gamma^2_{12}-\Gamma^2_{21})}{M}\,\nabla_{\boldsymbol\gamma}\boldsymbol\gamma
+k\,\tr(R)\ \boldsymbol\gamma.
\endgathered
\mytag{7.26}
$$
The derivatives $\nabla_{\boldsymbol\alpha}\boldsymbol\gamma$ and 
$\nabla_{\boldsymbol\gamma}\boldsymbol\gamma$ are taken from \mythetag{7.17}. 
The second order derivatives $\nabla_{\boldsymbol\alpha}\bigl(\nabla_{\boldsymbol
\gamma}\boldsymbol\gamma\bigr)$ and $\nabla_{\boldsymbol\gamma}
\bigl(\nabla_{\boldsymbol\alpha}\boldsymbol\gamma\bigr)$ are transformed as follows:
$$
\gather
\hskip -2em
\gathered
\nabla_{\boldsymbol\alpha}\bigl(\nabla_{\boldsymbol
\gamma}\boldsymbol\gamma\bigr)=
\nabla_{\boldsymbol\alpha}\bigl(\Gamma^1_{22}\,\boldsymbol\alpha+
\Gamma^2_{22}\,\boldsymbol\gamma
\bigr)=\nabla_{\boldsymbol\alpha}\bigl(\Gamma^1_{22}\bigr)
\,\boldsymbol\alpha\,+\\
\vspace{0.5ex}
+\,\Gamma^1_{22}\,\bigl(\Gamma^1_{11}\,\boldsymbol\alpha+
\Gamma^2_{11}\,\boldsymbol\gamma\bigr)+\nabla_{\boldsymbol\alpha}\bigl(
\Gamma^2_{22}\bigr)\,\boldsymbol\gamma
+\Gamma^2_{22}\,\bigl(\Gamma^1_{12}\,\boldsymbol\alpha
+\Gamma^2_{12}\,\boldsymbol\gamma\bigr),
\endgathered
\quad
\mytag{7.27}\\
\vspace{1ex}
\hskip -2em
\gathered
\nabla_{\boldsymbol\gamma}\bigl(\nabla_{\boldsymbol
\alpha}\boldsymbol\gamma\bigr)=\nabla_{\boldsymbol\gamma}\bigl(
\Gamma^1_{12}\,\boldsymbol\alpha+\Gamma^2_{12}\,\boldsymbol\gamma\bigr)
=\nabla_{\boldsymbol\gamma}\bigl(\Gamma^1_{12}\bigr)\,\,\boldsymbol\alpha\,+\\
\vspace{0.5ex}
+\,\Gamma^1_{12}\,\bigl(\Gamma^1_{21}\,\boldsymbol\alpha+
\Gamma^2_{21}\,\boldsymbol\gamma\bigr)
+\nabla_{\boldsymbol\gamma}\bigl(\Gamma^2_{12}\bigr)\,\boldsymbol\gamma
+\Gamma^2_{12}\,\bigl(\Gamma^1_{22}\,\boldsymbol\alpha
+\Gamma^2_{22}\,\boldsymbol\gamma\bigr).
\endgathered
\quad
\mytag{7.28}
\endgather
$$
Summarizing the formulas \mythetag{7.23}, \mythetag{7.24}, \mythetag{7.25}, 
\mythetag{7.26}, \mythetag{7.27}, and \mythetag{7.28}, we can formulate
the following lemma. 
\mylemma{7.1} The coefficients $F^1_1$, $F^2_1$, $F^1_2$, $F^2_2$ of the 
linear combinations \mythetag{7.22} are expressed through $M$, through
the pseudoscalar fields $\Gamma^1_{11}$, $\Gamma^2_{11}$, $\Gamma^1_{12}$, 
$\Gamma^2_{12}$, $\Gamma^1_{21}$, $\Gamma^2_{21}$, $\Gamma^1_{22}$, \
$\Gamma^2_{22}$ from \mythetag{7.17}, and through covariant derivatives
of them. 
\endproclaim
     Note that $\boldsymbol\alpha\nparallel\boldsymbol\gamma$ due to $M\neq 0$
in \mythetag{6.7}. Hence the values of the fields $\boldsymbol\alpha$ and 
$\boldsymbol\gamma$ constitute a basis at each point of $\Bbb R^2$ or a 
two-dimensional manifold. In this case the coefficients $F^1_1$, $F^2_1$, 
$F^1_2$, $F^2_2$ of the linear combinations \mythetag{7.22} constitute the
matrix of the linear operator $R$ in \mythetag{7.22}. It is well-known 
(see \mycite{12}) that the determinant of a linear operator does not depend 
on a basis where its matrix is calculated. Therefore we have the following
equality:
$$
\hskip -2em
\det(R)=\det\Vmatrix R^1_1 & R^1_2\\ 
\vspace{2ex}
R^2_1 & R^2_2\endVmatrix
=\det\Vmatrix F^1_1 & F^1_2\\ 
\vspace{2ex}
F^2_1 & F^2_2\endVmatrix.
\mytag{7.29}
$$
Applying \mythetag{7.29} along with Lemma~\mythelemma{7.1} to the formula
\mythetag{7.5}, we derive a theorem. 
\mytheorem{7.1} Within the intersection class\/ {\rm ShrID1\,$\cap$\,BgdET2}, 
i\.\,e\. if the conditions \mythetag{6.7} are fulfilled, Bagderina's pseudoscalar
field $j^{\,\sssize\text{Bgd}}_5$ from \mythetag{5.24} is expressed through
$M$, through $\Omega$, through the pseudoscalar fields $\Gamma^1_{11}$, 
$\Gamma^2_{11}$, $\Gamma^1_{12}$, $\Gamma^2_{12}$, $\Gamma^1_{21}$, 
$\Gamma^2_{21}$, $\Gamma^1_{22}$, $\Gamma^2_{22}$ from \mythetag{7.17}, 
and through covariant derivatives of them along the pseudovectorial fields 
$\boldsymbol\alpha$ and $\boldsymbol\gamma$. 
\endproclaim
\head
8. Comparison of invariants. 
\endhead
     Two scalar invariants in the first case of intermediate degeneration are 
very simple. They are given as two ratios of the pseudoscalar fields $M$, $N$ and 
$\Omega$:
$$
\xalignat 2
&\hskip -2em
I_1=\frac{M}{N^2},
&&I_2=\frac{\Omega^2}{N}
\mytag{8.1}
\endxalignat
$$
(see \thetag{6.8} and \thetag{6.19} in \mycite{3} or \thetag{5.1} in \mycite{4}).
Other scalar invariants are more complicated. The formulas \mythetag{7.17} and
the coefficients $\Gamma^1_{11}$, $\Gamma^2_{11}$, $\Gamma^1_{12}$, 
$\Gamma^2_{12}$, $\Gamma^1_{21}$, $\Gamma^2_{21}$, $\Gamma^1_{22}$, 
$\Gamma^2_{22}$ therein are used for introducing them (see \thetag{6.13} in 
\mycite{3} or \thetag{5.2} in \mycite{4}). Some of the coefficients 
in \mythetag{7.17} are identically zero:
$$
\xalignat 2
&\hskip -2em
\Gamma^2_{11}=0,
&&\Gamma^1_{21}=0.
\mytag{8.2}
\endxalignat
$$ 
Some others are expressed through the pseudoscalar field $N$:
$$
\xalignat 2
&\hskip -2em
\Gamma^1_{11}=-\frac{3}{5}\,N,
&&\Gamma^2_{21}=-\frac{3}{5}\,N.
\mytag{8.3}
\endxalignat
$$ 
And some of them are bound to the invariants $I_1$ and $I_2$ in \mythetag{8.1} and
to the pseudoscalar field $N$ by more complicated relationships
$$
\gather
\hskip -2em
I_1\,\Gamma^2_{12}=I_4\,N+\frac{3}{5}\,I_1\,N+2\,I_1^2\,N,
\mytag{8.4}\\
\vspace{2ex}
\hskip -2em
\aligned
 \bigl(I_1\,\Gamma^2_{22}\bigr)^4&+\bigl(I_7\,N^3\bigr)^2
 +\bigl(16\,I_2\,N^3\,{I_1}^4\bigr)^2=\\
 \vspace{1ex}
 &=32\,I_7\,N^6\,I_2\,{I_1}^4
 +2\,\bigl(I_7\,N^3+16\,I_2\,N^3\,{I_1}^4\bigr)\,
 \bigl(I_1\,\Gamma^2_{22}\bigr)^2
\endaligned
\mytag{8.5}
\endgather
$$
(see \thetag{6.22} and \thetag{6.23} in \mycite{3} or \thetag{5.6} and \thetag{5.7} 
in \mycite{4}). Unfortunately the formula \thetag{6.22} in \mycite{3} is mistyped
and then copied to \thetag{5.6} in \mycite{4}. The minus signs in the right hand side
of this formula should be altered for pluses. Here we present the correct formula
\mythetag{8.4}. As for the quantities $I_4$ and $I_7$ in \mythetag{8.5}, they are 
higher order invariants defined among others by the formulas 
$$
\xalignat 3
&\hskip -2em
I_4=\frac{\nabla_{\boldsymbol\alpha}I_1}{N},
&&I_5=\frac{\nabla_{\boldsymbol\alpha}I_2}{N},
&&I_6=\frac{\nabla_{\boldsymbol\alpha}I_3}{N},\\
\vspace{-1.5ex}
\mytag{8.6}\\
\vspace{-1.5ex}
&\hskip -2em
I_7=\frac{(\nabla_{\boldsymbol\gamma}I_1)^2}{N^3},
&&I_8=\frac{(\nabla_{\boldsymbol\gamma}I_2)^2}{N^3},
&&I_9=\frac{(\nabla_{\boldsymbol\gamma}I_3)^2}{N^3}.
\endxalignat
$$
Apart from \mythetag{8.2}, \mythetag{8.3}, \mythetag{8.4}, \mythetag{8.5}, we
have the relationship
$$
\hskip -2em
\Gamma^1_{12}=-\Gamma^2_{22}.
\mytag{8.7}
$$
(see \thetag{6.15} in \mycite{3}). Due to \mythetag{8.2}, \mythetag{8.3}, 
\mythetag{8.4}, \mythetag{8.5}, and \mythetag{8.7} the only nontrivial 
coefficient in \mythetag{7.17} is $\Gamma^1_{22}$. It is used in order to
produce the invariant $I_3$ in \mythetag{8.6}:
$$
\hskip -2em
I_3=\frac{\Gamma^1_{22}\,N^2}{M^2}
\mytag{8.8}
$$
(see \thetag{6.20} in \mycite{3}). The formula \thetag{5.5} for $I_3$ is
mistyped. The exponent of $M$ in the denominator is dropped. Here we
present the correct formula \mythetag{8.8}.\par
     In item 2 of her classification Theorem 2 Yu\.~Yu\.~Bagderina presents
two basic invariants. They are given by means of the formulas
$$
\pagebreak
\align
&\hskip -2em
I^{\,\sssize\text{Bgd}}_1=\frac{\Gamma^{\,\sssize\text{Bgd}}_0}
{\beta^{\,\sssize\text{Bgd}}_1\,(j^{\,\sssize\text{Bgd}}_0)^2},
\mytag{8.9}\\
\vspace{1ex}
&\hskip -2em
I^{\,\sssize\text{Bgd}}_2=\frac{5}{(j^{\,\sssize\text{Bgd}}_0)^2}
\bigl(2\,j^{\,\sssize\text{Bgd}}_1\,
j^{\,\sssize\text{Bgd}}_3+(j^{\,\sssize\text{Bgd}}_2
-j^{\,\sssize\text{Bgd}}_0/6)^2\bigr).
\mytag{8.10}
\endalign
$$
The first Bagderina's invariant is simple. Applying \mythetag{5.27},
\mythetag{5.28}, \mythetag{5.22}, we find 
$$
\hskip -2em
I^{\,\sssize\text{Bgd}}_1=\frac{N}{3\,\Omega^2}.
\mytag{8.11}
$$
\mylemma{8.1} Within the intersection class\/ {\rm ShrID1\,$\cap$\,BgdET2}, 
i\.\,e\. if the conditions \mythetag{6.7} are fulfilled, the first Bagderina's 
invariant $I^{\,\sssize\text{Bgd}}_1$ in \mythetag{8.9} is related to the 
invariant $I_2$ introduced in \mycite{3} according to the formula 
$$
\hskip -2em
I^{\,\sssize\text{Bgd}}_1=\frac{1}{3\,I_2}.
\mytag{8.12}
$$
\endproclaim
     Lemma~\mythelemma{8.1} and the formula \mythetag{8.12} in it are
immediate from \mythetag{8.1} and \mythetag{8.11}.\par
     The second Bagderina's invariant $I^{\,\sssize\text{Bgd}}_2$ in
\mythetag{8.10} is more complicated. Applying \mythetag{5.24} and 
\mythetag{5.22} to it, we can write the formula \mythetag{8.10} as follows:
$$
\hskip -2em
I^{\,\sssize\text{Bgd}}_2=\frac{j^{\,\sssize\text{Bgd}}_5}
{(j^{\,\sssize\text{Bgd}}_0)^2}=\frac{j^{\,\sssize\text{Bgd}}_5}
{9\,\Omega^2}.
\mytag{8.13}
$$
The numerator in \mythetag{8.13} is described by Theorem~\mythetheorem{7.1}. 
At this moment we know that the pseudoscalar fields $\Gamma^1_{11}$, 
$\Gamma^2_{11}$, $\Gamma^1_{12}$, $\Gamma^2_{12}$, $\Gamma^1_{21}$, 
$\Gamma^2_{21}$, $\Gamma^1_{22}$, $\Gamma^2_{22}$ from \mythetag{7.17}
are expressed through scalar invariants $I_1$, $I_2$, $I_3$, $I_4$, $I_7$
and the field $N$. Their covariant derivatives along $\boldsymbol\alpha$ and 
$\boldsymbol\gamma$ are expressed through the covariant derivatives of 
$I_1$, $I_2$, $I_3$, $I_4$, $I_7$ and through the covariant derivatives 
of $N$. \par
     Higher order covariant derivatives of  $I_1$, $I_2$, $I_3$ along 
$\boldsymbol\alpha$ and $\boldsymbol\gamma$ form higher order invariants
in \mythetag{8.6} and in the recurrent formulas 
$$
\xalignat 2
&\hskip -2em
I_{k+3}=\frac{\nabla_{\boldsymbol\alpha}I_k}{N},
&&I_{k+6}=\frac{\bigl(\nabla_{\boldsymbol\gamma}I_k\bigr)^2}{N^3}.
\mytag{8.14}
\endxalignat
$$
that should be applied in triples in some commonly negotiated order
(see \thetag{6.21} in \mycite{3}). Therefore covariant derivatives of 
$I_1$, $I_2$, $I_3$, $I_4$, $I_7$ of any order can be expressed back 
through the sequence of higher order scalar invariants, through $N$ and
through covariant derivatives of $N$.\par
     Let's consider the covariant derivatives along $\boldsymbol\alpha$ 
and $\boldsymbol\gamma$ for $N$ and for the other two fields $M$ and 
$\Omega$. In the case of $N$ we have
$$
\xalignat 2
&\hskip -2em
\nabla_{\boldsymbol\alpha}N=M,
&&\nabla_{\boldsymbol\gamma}N=-2\,M\,\Omega\,.
\mytag{8.15}
\endxalignat
$$
These formulas are derived by direct calculations in our special coordinates
introduced through Theorem~\mythetheorem{5.1}. In order to find covariant 
derivatives of $M$ we express it through $N$ and the scalar invariant 
$I_1$ by means of the first formula \mythetag{8.1}:
$$
\hskip -2em
M=I_1\,N^2. 
\mytag{8.16}
$$
Differentiating \mythetag{8.16}, we derive 
$$
\aligned
&\nabla_{\boldsymbol\alpha}M=(\nabla_{\boldsymbol\alpha}I_1)\,N^2+
2\,I_1\,N\,(\nabla_{\boldsymbol\alpha}N)=I_4\,N^3+2\,I_1\,N\,M,\\
&\nabla_{\boldsymbol\gamma}M=(\nabla_{\boldsymbol\gamma}I_1)\,N^2+
2\,I_1\,N\,(\nabla_{\boldsymbol\gamma}N)=
\sqrt{N^3\,I_7}\,N^2-4\,I_1\,N\,M\,\Omega\,.
\endaligned
\quad
\mytag{8.17}
$$
In the case of $\Omega$, we use the second formula \mythetag{8.1}. It
yields
$$
\hskip -2em
\Omega^2=I_2\,N.
\mytag{8.18}
$$
Differentiating \mythetag{8.18}, we derive the following formulas:
$$
\aligned
&\nabla_{\boldsymbol\alpha}\Omega=\frac{(\nabla_{\boldsymbol\alpha}I_2)
\,N+I_2\,(\nabla_{\boldsymbol\alpha}N)}{2\,\Omega}
=\frac{I_5\,N^2+I_2\,M}{2\,\Omega},\\
&\nabla_{\boldsymbol\gamma}\Omega=
\frac{(\nabla_{\boldsymbol\gamma}I_2)
\,N+I_2\,(\nabla_{\boldsymbol\gamma}N)}{2\,\Omega}
=\frac{\sqrt{N^3\,I_8\mathstrut}\ N-2\,I_2\,M\,\Omega}{2\,\Omega}. 
\endaligned
\quad
\mytag{8.19}
$$
Looking at \mythetag{8.15}, \mythetag{8.17}, and \mythetag{8.19}, we can 
formulate the following lemma.
\mylemma{8.2} Covariant derivatives of the pseudoscalar fields $M$, $N$,
and $\Omega$ along pseudovectorial fields $\boldsymbol\alpha$ and
$\boldsymbol\gamma$ are expressed through the scalar invariants $I_1$, $I_2$, 
$I_3$, $I_4$ etc in the recurrent sequence \mythetag{8.14} and through
these fields themselves. 
\endproclaim
    Note that Lemma~\mythelemma{8.2} applies not only to first order 
covariant derivatives, but to higher order derivatives as well, since
the formulas \mythetag{8.15}, \mythetag{8.17}, and \mythetag{8.19} can
be applied recursively. Combining this lemma with the above considerations
just after the formula \mythetag{8.13}, we derive a theorem.
\mytheorem{8.1}Within the intersection class\/ {\rm ShrID1\,$\cap$\,BgdET2}, 
i\.\,e\. if the conditions \mythetag{6.7} are fulfilled, the second Bagderina's 
invariant $I^{\,\sssize\text{Bgd}}_2$ in \mythetag{8.10} can be expressed
through $I_1$, $I_2$, $I_3$ and through higher order invariants $I_4$, $I_5$, 
$I_6$ etc in the recurrent sequence \mythetag{8.14}. 
\endproclaim
     The above proof of Theorem~\mythetheorem{8.1} is half constructive. One can
make it constructive by continuing the calculations ended with \mythetag{7.29},
though the resulting expression could be enormously large. 
\head
9. Conclusions.
\endhead
     Comparing the classification of the equations \mythetag{1.1} suggested 
by Yu\.~Yu\.~Bagderina in \mycite{5} with the previously known classification 
suggested in \mycite{3} we find that the case of intermediate degeneration 
from \mycite{3} does not completely coincide with the corresponding item 2 
of Bagderina's classification theorem in \mycite{5}. So, formally, 
Bagderina's classification is new. However, the case of intermediate 
degeneration from \mycite{3} has a substantial intersection with item 2 
in Bagderina's classification. Denoting the intersection class through
{\rm ShrID1\,$\cap$\,BgdET2}, we compared the two classifications within this
intersection class. As a result we have found that most basic structures 
and basic formulas from Bagderina's paper \mycite{5} do coincide or are 
very closely related to those in \mycite{3} and \mycite{4}, though they are 
given in different notations (see Lemma~\mythelemma{5.1}, Lemma~\mythelemma{5.2},
Lemma~\mythelemma{5.3}, Lemma~\mythelemma{5.4}, Lemma~\mythelemma{5.5},
Lemma~\mythelemma{5.7}, Lemma~\mythelemma{6.1}, Lemma~\mythelemma{6.2}, 
and Lemma~\mythelemma{8.1}).\par
     In item 2 of her Teorem 2 in \mycite{5} Yu\.~Yu\.~Bagderina presents 
two basic invariants $I^{\,\sssize\text{Bgd}}_1$ and $I^{\,\sssize\text{Bgd}}_2$, 
while in \mycite{3} three basic invariants $I_1$, $I_2$, and $I_3$ were 
presented. Yu\.~Yu\.~Bagderina claims that her two invariants are sufficient
for expressing all of the invariants, including $I_1$, $I_2$, $I_3$, through 
them and through their invariant derivatives. However, in \mycite{5} there
are no explicit formulas expressing $I_1$, $I_2$, $I_3$ through 
$I^{\,\sssize\text{Bgd}}_1$ and $I^{\,\sssize\text{Bgd}}_2$. Some formulas
of this sort are given in \mycite{13}, but again for the case $\Omega=0$, 
which is outside our present intersection class {\rm ShrID1\,$\cap$\,BgdET2}.
\par
     In the present paper we solve the basic invariants problem from our
side by showing that both Bagderina's invariants $I^{\,\sssize\text{Bgd}}_1$ 
and $I^{\,\sssize\text{Bgd}}_2$ can be expressed through the invariants $I_1$, 
$I_2$, $I_3$ and through proper invariant derivatives of them (see 
Lemma~\mythelemma{8.1} and Theorem~\mythetheorem{8.1} above). It would be best 
if Yu\.~Yu\.~Bagderina presents some explicit formulas or an algorithm for 
expressing $I_1$, $I_2$, $I_3$ through her invariants $I^{\,\sssize\text{Bgd}}_1$ 
and $I^{\,\sssize\text{Bgd}}_2$ in the intersection class 
{\rm ShrID1\,$\cap$\,BgdET2}. Otherwise her claim that her invariants are 
basic is open to question.
\par

\enddocument
\end